\documentclass[aps,preprint,a4paper,nofootinbib,floatfix]{revtex4-1}
\usepackage{graphicx}
\usepackage{amsmath,amssymb,amsfonts,amsthm,amscd,bm}
\usepackage{booktabs}
\usepackage{listings}
\usepackage{xcolor,colortbl}
\usepackage{float}

\def\bar{\overline}

\def\*{\star}
\def\[{\left[}
\def\]{\right]}
\def\({\left(}

\def\){\right)}

\def\frac#1#2{{#1 \over #2}}
\def\inv#1{{1 \over #1}}

\def\2pi{\hbox{$2\pi i$}}

\def\dsl{\raise.15ex\hbox{/}\kern-.57em\partial}
\def\Dsl{\,\raise.15ex\hbox{/}\mkern-.13.5mu D}

\def\CD{{\cal D}}

	\def\CN{{\cal N}}	
\def\CP{{\cal P}}		
\def\CS{{\cal S}}

\def\barray{\begin{eqnarray}}
\def\earray{\end{eqnarray}}
\def\beq{\begin{equation}}
\def\eeq{\end{equation}}

\def\smallhalf{\tfrac{1}{2}}
\def\cs{b}

\def\Sexact{X}
\def\Eu{{\rm \CP}}

\def\Ei{{\rm Ei}}

\makeatletter
\g@addto@macro\bfseries{\boldmath}
\makeatother

\begin{document}

\title{On the validity of the Euler product  inside the critical strip}

\author{Guilherme Fran\c ca\footnote{guifranca@gmail.com} and  
Andr\' e  LeClair\footnote{andre.leclair@gmail.com}}

\affiliation{Cornell University,  Physics Department, Ithaca, NY 14850} 


\begin{abstract}
The Euler product formula relates Dirichlet $L(s,\chi)$  functions to
an infinite product over primes,  and is known to be  valid for $\Re (s) >1$,
where it converges absolutely. 
We provide  arguments that 
the formula is actually  valid for $\Re (s) > 1/2$
 in a specific sense.
Namely,  the logarithm of the Euler product,   
although formally divergent, is 
meaningful because it is Ces\`aro summable,  and   
its   Ces\`aro  average  converges to $\log L (s,\chi)$. 
Our argument relies on
the prime number theorem,   an Abel transform,  
and  a central limit theorem for the 
Random Walk of the Primes, the series
$\sum_{n=1}^N  \cos\(t\log p_n\)$, and its generalization to other
Dirichlet $L$-functions.        
The significance of ${\Re(s) > 1/2}$ arises from the $\sqrt{N}$ growth
of this series, since it satisfies a central limit theorem.
$L$-functions based on principal Dirichlet characters,
such as the Riemann $\zeta$-function,  are exceptional
due to the pole at $s=1$,  and require $\Im (s) \neq 0$ 
and a truncation of the Euler product.  
Compelling numerical evidence of this surprising result is
presented, and some of its  consequences are discussed.  
\end{abstract}

\maketitle

\section{Introduction and Summary}

The Riemann $\zeta$-function was originally defined by the series
\beq
\label{zetaseries}
\zeta (s)  =   \sum_{n=1}^\infty \dfrac{1}{n^s} 
\eeq
where $s=\sigma + it$ is a complex number.  
This series converges absolutely for $\Re (s) > 1$.   It can be analytically 
continued to the entire complex plane by extending an integral representation
valid for $\Re (s) > 1$, except for the simple pole at $s=1$.   
Using only the unique prime factorization theorem,  one can derive the 
Euler product formula,   
which is the equality  
\beq
\label{Eu}
\zeta(s) =  \Eu(s) \equiv
\prod_{n=1}^\infty  \(  1 - \dfrac{1}{p_n^{\,s}} \)^{-1}
\eeq
and $p_n$ is the $n$th prime number.  
It is this formula which  is the key to Riemann's result \cite{Riemann}
that relates the 
distribution of primes, namely 
the prime number counting function $\pi (x)$,  to a series involving an 
infinite sum over zeros $\rho$ of
the $\zeta$-function inside the critical strip $0<\Re (s) <1$.  
Henceforth, it is implicit that 
$\zeta (s)$  inside the strip is defined by analytic continuation 
in the standard way.

The product \eqref{Eu} also  converges absolutely  only for $\Re (s) > 1$, 
and is known not to even conditionally converge for $\Re (s) \le 1$.  
However,
in their own studies, Berry and Keating used the Euler product 
inside the strip \cite{Berry}.  
Away from the real line, in some region of the critical strip,   
the phases $e^{i\, t\log p_n}$  can be such that 
one can make sense of the infinite product,  and this is the main idea  
studied in this article.     
Just to illustrate, 
if one introduces alternating signs into the series \eqref{zetaseries},  it
converges for $\Re (s) > 0$ and is the Dirichlet $\eta$-function. 
In this case, it converges by the simple alternating series
test.   
But if the signs were not strictly alternating,   convergence 
would be difficult to disprove or prove.   

The Riemann $\zeta$-function is the simplest,  trivial 
example of $L$-functions based on Dirichlet 
characters $\chi$,   and it will be important to consider the
full class of such $L$-functions. 
Here the Euler product formula takes the form
\beq
\label{EuL}
L(s,\chi) = \sum_{n=1}^\infty  
\dfrac{\chi(n)}{n^s} = \CP (s,\chi) \equiv  \prod_{n=1}^{\infty} 
\(1 - \dfrac{\chi (p_n)}{p_n^{\,s}}\)^{-1} .
\eeq
  
In this article we provide strong arguments 
(although not a strict mathematical proof) 
that the Euler product formula 
is  valid in a concrete 
sense for $\Re (s) > 1/2$,  i.e. 
on the right-half of the critical strip.  
Our  argument invokes 
the prime number theorem,
an  Abel transformation (summation by parts),  and
the \emph{Central Limit Theorem} (CLT) for the particular 
series \eqref{BNdef} below.       
The most important ingredient is the CLT,
and the significance of  $\Re(s) = 1/2$ comes from 
the $\sqrt{N}$ growth of the series \eqref{BNdef}.    
We emphasize that \emph{we do not introduce any probabilistic 
aspect to the original problem},  and this work is not in the realm of
so-called probabilistic number theory.     The CLT is
only invoked as a tool in order  to establish this $\sqrt{N}$ growth;   
the original series 
is a unique deterministic  member of the ensembles of the 
CLT\footnote{
An amusing quote
of Poincar\'e is relevant:  \emph{``there must be something 
mysterious about the normal law since mathematicians think 
it is a law of nature whereas physicists are convinced that 
it is a mathematical theorem.''} 
In the present work there is no data from nature,  and the CLT is
indeed a mathematical theorem.}. We will also provide
compelling  numerical evidence of this surprising result.  

The Euler product \eqref{Eu} was previously studied inside
the critical strip by Gonek, Hughes, and Keating \cite{Gonek}.
They proved that $\zeta(s)$ can be well-approximated
by a truncated hybrid Euler-Hadamard product for $\sigma > 0$; 
$\zeta(s) = \CP_X(s)  \, Z_X(s) \big(1+o(1)\big) $   where
$\log \CP_X(s) = \sum_{n \le X} \Lambda(n) / (n^s \log n)$ with
$X < t^{1-\epsilon}$ and $|s-1| > 10^{-1}$. This sum
is related to a truncated Euler product \eqref{Eu} through the relation
$\tfrac{d}{ds}\log \zeta(s) = \sum_{n\ge1}\Lambda(n) n^{-s}$, valid
for $\sigma > 1$, but it was extended into the strip by 
Titchmarsh \cite{Titchmarsh}. 
The function $Z_X(s)$ is a Hadamard product that depends 
on all non-trivial zeros of $\zeta(s)$. This result was also extended
to Dirichlet $L$-functions \cite{Keating1,Keating2,Keating3}.
Gonek \cite{Gonek1} improved on the result \cite{Gonek}
by introducing a smoothed
sum, $\log \CP_X(s) = \sum_{n \le X^2} \Lambda_X(n)/(n^s \log n)$, 
and also
proved that {\it  assuming the Riemann Hypothesis (RH)}, there is no
contribution from the Hadamard product for $\sigma > 1/2$, 
i.e. $Z_X(s)\to 1$,  and then
$\zeta(s) = \CP_X(s)\big( 1+O(\log^{(1-C)/2} t) \big)$ 
with $2 \le X \le t^2$, $|s-1| \ge 1/10$ and
$1/2+C \log_2 (2|t|+4) / \log X \le \sigma \le 1$, for some constant $C > 1$. 
Thus, under the RH, a short truncation of the Euler product \eqref{Eu}
approximates $\zeta(s)$ into the right-half part of the critical strip, but not
too close to the critical line.
In this paper, we are going to analyze the Euler product from a different
perspective, namely, we are going to analyze the convergence behaviour
of \eqref{Eu} directly, using
properties of the primes.
In other words,  our starting point
is not the formula $\zeta(s) \approx \CP_X(s)  Z_X(s)$ since we 
assume nothing about the zeros, and we have no analog of $Z_X(s)$. 
Furthermore,  we will not assume the RH.   

It is important to mention that partial Euler 
products \emph{on the critical line}
of more general $L$-functions 
was also considered by Conrad \cite{Conrad}, 
and in the case of nontrivial Dirichlet $L$-functions
an interesting theorem equivalent to the
RH was demonstrated, related to the product at ${s=1/2}$. 
The behaviour of Euler products of Dirichlet $L$-functions
on the critical line was also studied in \cite{Kimura}.

In this paper we are interested in determining what is the largest
region in the critical strip that the Euler product can be meaningful in
its most basic sense.
We first prove the following  in 
the next section. 
Consider the series
\beq
\label{BNdef}  
B_N(t,\chi) = \sum_{n=1}^N  \cos ( t \log p_n - \theta_n ),  \qquad
\theta_n = -i \log \chi (p_n),
\eeq
which we will refer to as the 
\emph{Random Walk of the Primes} (RWP), even though 
it is  a completely deterministic series.
If $B_N$ grows as $\sqrt{N}$,  then  $\log \Eu(s,\chi)$  
converges for $\sigma > 1/2$ if one takes its Ces\`aro average.
This Ces\`aro average is not probabilistic,  but rather is a 
smoothing procedure.

This $\sqrt{N}$ growth is robust and universal in  statistics 
and statistical physics. For instance, diffusion grows as the square-root 
of time. For  a random variable with standard deviation $\sigma$, 
the relative uncertainty goes   as $\sigma / N \sim 1/\sqrt{N}$ and 
thus becomes small for large $N$.   
This is  a consequence of the CLT. 
We will establish  that  the CLT 
applies to \eqref{BNdef}.      
The specialty of $\Re(s) =1/2$  is due to this square root. 
The beauty of this argument is that it does not 
rely on any details of the primes,  on the contrary,  it depends on 
their \emph{multiplicative independence}, which is reflected in their 
pseudo-random  behaviour \cite{Tao}. 
This is analogous to the fact that one does not need to know the exact
positions and velocities of $N\sim 10^{24}$ 
molecules in a gas  to predict its pressure.   The situation for the RWP is
even better compared to this,  
since the list of primes is very long, 
especially towards the end.   

The Euler product  in \eqref{Eu} does not converge  
inside the critical strip  in the conventional sense,  
since the domains of convergence of Dirichlet series are always 
half-planes, and due to the pole of $\zeta(s)$ at $s=1$, it  implies that
it can only converge for $\sigma > 1$.
However, some divergent series  are  still 
meaningful \cite{Hardy}.   A formally  divergent series can
still be summable,   if the 
divergence simply amounts to 
fluctuations  around a meaningful central value.  
This is referred to as Ces\`aro summability,  which means that
its average converges.  
More precisely,   
in the sequel we will provide arguments for  the following equality
\beq
\label{EPFave}
\log L  (s, \chi )  =  \langle \log \CP(s,\chi) \rangle \qquad
\mbox{for $\Re (s) > 1/2$,}
\eeq
where $\langle \log \CP(s,\chi) \rangle$ denotes its Ces\`aro average, 
and $L (s, \chi)$ is the standard analytic continuation of the series
\eqref{EuL} into the critical strip.
As we will explain,  for $\zeta(s)$ and  other principal $L(s,\chi)$,  
the above equation needs to be refined by introducing a cut-off 
$N_c$, truncating the product,  and the $=$ sign should
be replaced by $\approx$.   
The cut-off can only be taken to infinity in the limit $t \to \infty$.
This is due to the pole at $s=1$,  
which does not exist for non-principal Dirichlet $L$-functions.

By EPF  let us  refer to the Euler product formula \eqref{Eu},
or \eqref{EuL}, for $\Re(s) > 1/2$.
If it is indeed valid, even in the average sense described here,
there are many consequences.
There is one which is immediate.
It is well-known that  Euler product formula 
implies that $L(s,\chi)$ has no zeros with $\Re (s) >1$.  
The same argument applies to \eqref{EPFave}:    
if $\langle \log \CP(s,\chi) \rangle$ is finite,   then
$\log L (s, \chi)$  is never infinite.   On the other hand,  a
zero $\rho$ of $L(s,\chi)$
implies $\log L(\rho, \chi) = -\infty$,   thus 
there  are no zeros with $\Re (s) >1/2$.    
Incidentally, the EPF also gives a new  proof of
the  prime number theorem,
which is equivalent to the fact that
there are no zeros of $\zeta (s)$  with $\Re (s) =1$.     
In fact,   nothing very special happens while crossing 
$\Re (s) =1$;   in contrast the behavior changes
dramatically at $\Re (s) =1/2$.    
The $\zeta$-function satisfies the functional equation
\beq
\label{chi}
\chi (s) = \chi(1-s) = \pi^{-s/2} \Gamma (s/2)  \zeta (s), 
\eeq
and Dirichlet $L$-functions satisfy a similar equation.  
This then shows there are also no non-trivial zeros with $\Re (s) < 1/2$.
Thus the EPF combined with the functional equation constrains all
non-trivial  zeros to be on the 
critical line $\Re(s)=1/2$, 
which is of course the 
\emph{Generalized Riemann Hypothesis} (GRH).     
    
We are proposing that it is ultimately the multiplicative independence of the
primes,  together with the strongly multiplicative property of  
Dirichlet characters, which makes the series \eqref{BNdef} behave like 
a sum of
independent random variables, that underlies the validity of the GRH.
Other  consequences will be discussed in the last section of this article.  
For one,  it provides further validation of the transcendental equations for
individual zeros derived in \cite{RHLeclair,FrancaLeclair}. 
It also leads to  a formula  that relates Riemann zeros to an 
infinite sum over primes, 
which is a  kind of inverse of Riemann's result that relates primes to sums
over zeros.

We organize our work as follows. In Section~\ref{sec:convergence},
we give a criterion for the convergence of the average of the 
Euler product \eqref{EuL}.
In Section~\ref{sec:sqrtN}, we show
that \eqref{BNdef} obeys a CLT, which implies the bound 
$B_N = O(\sqrt{N})$. We also discuss the 
difference in  behaviour of \eqref{BNdef}  for 
non-principal verses  principal characters.  For the later, we need
to introduce a cut-off, truncating the sum. This is fundamentally
related to the pole of $L(s,\chi) $ at $s=1$. 
The same applies to the $\zeta$-function.
In Section~\ref{sec:transcendental}, we discuss
some consequences of the Euler product formula inside
the right-half part of the critical strip. More precisely, 
the implications to the transcendental equations for the $n$th non-trivial
zero, proposed in \cite{FrancaLeclair}. 
This is due essentially to the behaviour of 
the argument of $L(s,\chi)$ on the critical line, which
can then be described through the Euler product. We 
present our final considerations 
in Section~\ref{sec:conclusion}. 
The Appendix~\ref{sec:numerical} contains
numerical results, validating our statements.

\section{A criterion for finiteness  of the Euler Product}
\label{sec:convergence}

The product in \eqref{EuL} converges if the following sum converges
\beq
\label{logsum}  
\log \Eu(s,\chi) =  
- \sum_{n=1}^\infty \log \( 1- \dfrac{\chi (p_n) }{p_n^{\,s}} \)  = 
\sum_{n=1}^{\infty} \(  \dfrac{\chi(p_n)}{p_n^{\,s}} +   
\dfrac{\chi(p_n)^2 }{2 p_n^{\,2s}} + \dotsm  \)  .
\eeq
The second term  and higher in \eqref{logsum} converge 
absolutely for $\sigma>1/2$.  
Thus convergence of the Euler product depends on the first term, 
i.e. on the series
\beq
\label{Series1}
\Sexact (s,\chi)  =  \lim_{N\to \infty}  \Sexact_N (s,\chi)   =  
\lim_{N\to \infty}   \sum_{n=1}^N   
\dfrac{\chi (p_n) }{(p_n)^{s} }.
\eeq
Chernoff \cite{Chernoff} considered the above series  
for the trivial character $\chi =1$, with $p_n$ 
replaced by $n\log n$, and showed it could be analytically
continued for $\sigma >0$;  therefore, the hypothetical zeta function based
on this product has no zeros in the entire critical strip.   
As we will see,  it is important not to do this in the 
phase.  For $\chi =1$,  the series \eqref{Series1} is known as the 
\emph{Prime Zeta Function}, and has a rich pole structure. 
It can be analytically continued to $\sigma > 0$, except
for poles on the real line $0<s \le 1$, and points corresponding to 
Riemann zeros. The imaginary line $\sigma = 0$ is a natural boundary
of the function.

Already we saw  the role of $\sigma = 1/2$,  but for elementary 
reasons that are clearly not enough for our purposes.  
In the trivial case $\chi=1$ corresponding to the $\zeta$-function, 
the series \eqref{Series1} only  
converges absolutely for $\sigma >1$.   
Actually, it  also fails the Dirichlet test of convergence 
since $|\sum_n e^{-i t \log p_n }|$  is unbounded;   
if it were bounded then the series would
converge for all $\sigma >0$,  which is certainly not the case,  otherwise 
this would rule out  the known infinite number of Riemann zeros on the 
critical line, and also the pole at $s=1$. 
Thus \eqref{Series1} fails the  simplest  convergence tests.

Let us consider convergence 
of the real and imaginary parts of \eqref{Series1} separately.
Let $\CS(s,\chi)$ denote the real part of \eqref{Series1}. Then
we have
\beq
\label{NewS} 
\CS(s,\chi) =   \lim_{N\to \infty}   
\CS_N(s,\chi)   =  \lim_{N\to \infty} \sum_{n=1}^N a_n b_n 
\eeq
where
\beq
\label{anbn}
a_n =   p_n^{\,-\sigma} , \qquad
b_n =   \cos ( \lambda_n ),  \qquad 
\lambda_n = t \log p_n - \theta_n, \qquad 
\theta_n = - i \log \chi (p_n). 
\eeq
The characters $\chi$ are all either a phase,  or zero,  
thus $\lambda_n$ is real.   
It is implicit that terms corresponding to $\chi(p_n) =0$
are omitted in the above sum, since they do not contribute to \eqref{Series1}.
Analogous arguments  apply to the imaginary part 
of \eqref{Series1},  with  
$b_n = -\sin(\lambda_n)$.
As stated above, for the $\zeta$-function with $\chi=1$,  
when $t=0$, 
then \eqref{NewS} converges only if $\sigma >1$.
However, in the general case, the oscillations  of $\cs_n$ can conspire 
to make  the series \eqref{NewS} converge  for  $\sigma \le 1$.
The simplest example to 
illustrate this is to replace $\cs_n$ by $(-1)^n$. 
In this case the alternating sign test
shows that the series converges for $\sigma > 0$.      
For our series \eqref{NewS},  
the signs of $\cs_n$ can be both positive or negative,  
but they do not strictly alternate. 
Rather, the situation here is between the two extremes of strictly 
alternating signs 
verses all positive signs,    
which suggests that \eqref{NewS} may converge for 
$\sigma > \sigma_c$ for some $\sigma_c \in (0,\,1]$.

Through an Abel transformation, the partial sum 
in \eqref{NewS} can be rewritten as 
\beq
\CS_{N} =  a_{N} B_{N}   - \sum_{n=1}^{N-1}  B_n (a_{n+1} - a_n)
\qquad \mbox{where} \qquad
B_n = \sum_{k=1}^n  b_k .
\label{BN}
\eeq
This implies 
\beq
\label{sumparts2}
|\CS_{N}|  <   | a_{N}| | B_{N}|   +  
\sum_{n=1}^{N-1}  |B_n|  |a_{n+1} - a_n|.
\eeq
Now, we  have that
\beq
\label{angap}
|a_{n+1} - a_n| = 
\left| \dfrac{1}{(p_{n+1})^\sigma} - \dfrac{1}{(p_n)^\sigma}\right| < 
\sigma \frac{g_n}{p_n (p_{n+1})^\sigma}
\eeq
where $g_n = p_{n+1} - p_n$ is the gap between consecutive
primes.

Let us now use two aspects of the prime number theorem. 
The first is that $p_n >  n \log n$. 
The second is that on average $\langle p_n  \rangle  = n \log n$, 
thus the  
\emph{average gap} is $\langle g_n \rangle = \log n$.
Thus, on average 
\beq
\label{AbsSN} 
\langle|\CS_N| \rangle < \frac{|B_N|}{(N \log N)^\sigma} + \sum_{n=1}^{N-1}
\dfrac{\sigma |B_n|}{n (n\log n)^\sigma}.
\eeq

For the moment,  let us assume that $|B_N|$ grows as $\sqrt{N}$,
i.e. $B_N = O(\sqrt{N})$,  which will be justified in the next section.      
Then,  
the RHS of the above equation behaves like   $\sum_n 1/n^{\sigma + 1/2}$,
which implies that the average  $\langle 
\CS_N \rangle$  converges  for $\sigma> 1/2$.   
There are many ways to define   $\langle \CS_N \rangle$,  
since the average is really just a smoothing procedure,  
but if the central value is  meaningful,  they  should agree.   
The simplest way is   to replace  
it   by its arithmetic average 
\beq
\label{aveS}
\langle \CS_N  \rangle =  \dfrac{1}{N} \sum_{n=1}^N  \CS_n .
\eeq 
One can easily show that in the limit of large $N$, then $\langle \CS_{N+1}
\rangle = \langle  \CS_N \rangle $.   
Convergence in this average sense  is referred to as
Ces\`aro summability.   
Henceforth,  we will refer to 
this  convergence  in an average sense as Ces\`aro-convergence,
and unless otherwise stated, simply convergence for short.

\section{$\sqrt{N}$ growth of $B_N$ from a central limit theorem}
\label{sec:sqrtN}

In the last section,  we showed that if   
$B_N$ grows as $\sqrt{N}$,  then
the Euler product  \eqref{EuL} Ces\`aro converges  on the right-half of the 
critical strip. 
In this section we  
prove that, for \emph{non-principal} characters,  
$\lim_{N\to \infty} B_N/\sqrt{N}$ is finite by 
using a version of the central limit theorem.
As we will see,  for \emph{principal} characters this statement needs to
be refined.  

\subsection{The general result}

The simplest, and original, version of the CLT is for independent and  
identically distributed (iid) random variables.
Let us recall the statement of the theorem in this case. Consider
\beq
\label{Rseries}
R_N = \sum_{n=1}^N  r_n 
\eeq
where $r_n$ are iid random variables with 
\emph{zero mean} and \emph{finite variance}.
For example,  if  $r_n=\pm 1/\sqrt{2}$, 
then this series is the standard 
random walk in one dimension.  
Below,  we will  consider $r_n$  as a real  random variable  
uniformly distributed on the interval $[-1,1]$.
In either case the distribution of $R_N/\sqrt{N}$ 
approaches a gaussian (normal) distribution at large $N$,
with zero mean and variance $1/2$, namely
\beq
\label{prob}
\lim_{N\to\infty} \dfrac{R_N}{\sqrt{N}} \to \CN_{0,1/2}
\eeq
where $\CN_{0,1/2}$ is a normal random variable determined
by a gaussian density $\tfrac{1}{\sqrt{\pi}}e^{-u^2}$.
The CLT guarantees that in the limit of large $N$, 
then $R_N /\sqrt{N}$ is finite for any member of the ensembles.      

What is important for our purposes is that certain 
trigonometric series are known to behave as  
iid random variables and thus satisfy a CLT. 
Consider the series
\beq
\label{CN}
C_N (u) = \sum_{n=1}^N \cos (u \, \lambda_n ) 
\eeq
where $u$ is a uniformly distributed real variable on the interval 
$[0, 2\pi]$. A well known example is the lacunary trigonometric series 
\cite{Salem}, 
where $\lambda_n$ are integers with gaps that grow fast enough,
namely they satisfy the Hadamard gap condition, 
$\lambda_{n+1}/\lambda_n  >  q >1$ for all $n$.   
For example, $\lambda_n= 2^n$ satisfies this condition.   
Clearly, for a fixed $u$, the terms in the series \eqref{CN}  are 
not  iid random variables 
since the $\lambda_n$'s are deterministic and highly correlated,
nevertheless the CLT is still valid. Although the theorem originally 
assumed that  $\lambda_n$ is an integer,  
it was later shown that this
is an unnecessary restriction \cite{Salem2}.  
Our series $B_N$ is 
equal to $C_N (u=1)$ with $\lambda_n$ given in \eqref{anbn}.   
Unfortunately, one cannot apply the theorems  for lacunary
trigonometric series since these $\lambda_n$'s do not satisfy the
Hadamard gap condition.
 
Let us first present some heuristic arguments before stating 
a precise result.   
The primes are deterministic,  
nevertheless, it is generally 
accepted that they behave pseudo-randomly\footnote{
{\it ``God may not play dice with the universe, but something strange is going 
on with the prime numbers''}. This is a misattributed quotation to
P. Erd\H os, one of the pioneers in applying probabilistic methods
to number theory, but actually it seems to be a comment from 
Carl Pomerance in a talk about the Erd\H os-Kac theorem, in response to
Einstein's famous assertion about quantum mechanics.} \cite{Tao}.
If the primes were truly random,   
then since $-1 \leq \cos ( \lambda_n ) \leq 1$,
the series $B_N$ should behave 
like $R_N$, with $r_n$ uniformly distributed
on the interval $[-1,1]$.   As we will show,  this heuristic argument 
leads to the correct result.   
However,  the problem with the above  argument is that 
pseudo-randomness of the primes is a somewhat  vague concept, and 
difficult to quantify.   

Fortunately, one can prove the desired 
result using only the \emph{multiplicative  independence}
of the primes, and the strongly multiplicative property of Dirichlet
characters, which is essential to derive the Euler product formula.
It is known that the CLT applies to the series 
\eqref{CN}  if the $\lambda_n$'s  are  linearly independent 
over the integers;  see for instance 
\cite[pp. 47]{Kac} and \cite[pp. 35]{Billingsley2}.
It is easy to show that our $\lambda_n$ have this property. 
For any integer $I>1$,  from the unique prime factorization theorem
one has
\beq
\label{IDir}
I^{\,t}  \cdot \chi(I)^{i}  =  
\prod\nolimits_k (p_k)^{t\, n_k} \,  \chi (p_k)^{i n_k} 
\eeq
where $n_k$ are integers.  
Taking the logarithm one finds
\beq
\label{indepDir}
\sum_k  n_k \( t\, \log p_k  + i \log \chi (p_k) \) =
\sum_k   n_k \, \lambda_k  \neq 0
\eeq
where it is implicit that terms with $\chi(p_k) =0$ are dropped from
the sum.  
Therefore, the $\lambda_n$'s in \eqref{anbn} are linearly independent. 
Thus the series 
\beq
\label{Bprime}
B_N (u; t, \chi) =  \sum_{n=1}^N  \cos \[ u(t \log p_n - \theta_n )   \]
\eeq
satisfies the CLT.
The original series \eqref{BN} of the last section, see also \eqref{BNdef},
corresponds to $u=1$. It is useful to introduce the 
additional variable $u \in [0,2\pi]$  since it allows  us to 
study the distribution of $B_N (u)$ on a 
given interval.  
We emphasize that this is simply a useful device and it does not
introduce any additional probabilistic aspect to the original series 
$B_N$, equation \eqref{BNdef},  which is completely deterministic.
The  CLT for \eqref{Bprime} guarantees that 
$\lim_{N\to \infty}
B_N / \sqrt{N}$  is finite for any $u\neq 0$.
Thus the Euler product for $L(s, \chi)$  Ces\`aro 
converges for $\Re (s) >1/2$, according to the discussion
of the previous section.

There is a very  important  difference  
between $L(s, \chi)$ with \emph{principal}  verses \emph{non-principal} 
characters.  The characters can be denoted as $\chi_{k,j}$, 
where $k$ is the modulus,  and $j = 1, \ldots, \varphi (k)$
where $\varphi (k)$ is the Euler totient, and equals the number of
distinct characters of modulus $k$.    For each $k$ there is only 
one principal character,  denoted by $\chi_{k,1}$,  which  is defined as
$\chi_{k,1} (n) = 1$ if $k$ and $n$ are coprime, and $\chi_{k,1} (n) = 0$ 
otherwise.
The $\zeta$-function corresponds to the Dirichlet $L$-function 
for the trivial principal character
of modulus $k=1$, where $\chi_{1,1} (n) =1$ for every $n$.
In fact, the non-trivial zeros of all $L$-functions based on
principal characters are the same as for $\zeta(s)$.   
For a principal character,  the terms which contribute to the
sum \eqref{Series1} are $\chi (p_k) = 1$, implying 
that $\lambda_k = 0$, unless $t\neq 0$.
Thus for the principal characters,   there is no CLT 
whatsoever for $t=0$.    
This implies that for $t=0$ we have the growth
$B_N = O(N)$, and the series  actually \eqref{NewS} diverges.

The origin of the latter divergence 
is of course  the existence of a pole
at $s=1$ for principal Dirichlet $L$-functions.
On the other hand,  the vast majority of 
Dirichlet $L$-functions are non-principle,  since
$\varphi (n)$ increases with $n$,  and $\lambda_k \ne 0$ 
even when $t = 0$. In the case $t=0$, then 
$\lambda_k = i \log \chi (p_k)$, 
and they are still linearly independent by \eqref{indepDir}.   
Therefore, for all Dirichlet $L$-functions, except for those based on 
principal characters, 
the Euler product Ces\`aro-converges for $\Re (s) >1/2$, 
\emph{including the real line}.
This is consistent, and  in fact predicts,  that  unlike the 
$\zeta$-function,  these $L$-functions have no poles   on the 
real line $\Re (s) >1/2$,  which is known to be the case.

\subsection{The natural cut-off for principle characters}
\label{sec:cutoff}

Clearly, the CLT holds for non-principal characters as discussed
in the previous section. In this case we can take
$N$ arbitrarily large, regardless of the value of $t$. 
This can be intuitively seen as
follows. Consider the phase difference between two consecutive
waves in \eqref{Bprime} (with $u=1$),
\beq\label{deltaphiD}
\Delta \phi = t \log p_{n+1} - t \log p_n - (\theta_{n+1} - \theta_n)
\approx t \dfrac{g_n}{p_n} + \( \theta_{n} - \theta_{n+1} \).
\eeq
For a fixed $t$ and large $N$, $t g_N/p_N \to 0$. However, for non-principal
characters, there is always a difference of phase due to the second
term in \eqref{deltaphiD} that allows cancellations between different
waves. The waves are not going to suffer a totally constructive
interference, and this is the reason for a
growth of $B_N = O(\sqrt{N})$, instead of $B_N = O(N)$.
Note that we can even set $t=0$ without problem.

The situation is more complicated for the $\zeta$-function, and also for
principal Dirichlet \mbox{$L$-functions}. In these two cases
we have $\chi(p_n) = 1$ for all $n$, and the second term
$(\theta_{n+1}-\theta_n)$ is not
present in \eqref{deltaphiD}.
Therefore, there is a subtlety  in applying the CLT to the series $B_N$ when 
there is no contribution from the characters,  namely the series
\beq\label{BN2}
B_N(t)  =  \sum_{n=1}^N  \cos ( t \log p_n).
\eeq
For unknown reasons, there is no mention of
such a subtlety in the original works \cite{Kac,Billingsley}.
As we now discuss,  
it is possible to still have $B_N = O(\sqrt{N})$  valid  
below a cut-off $N_c$, that depends on $t$.   
Strictly speaking it is not a CLT since $N$ cannot be taken 
freely to infinity.
The phase difference between consecutive waves in \eqref{BN2} is
\beq\label{deltaphi}
\Delta \phi = t \log p_{n+1} - t \log p_n \approx t \dfrac{g_n}{p_n}.
\eeq
When $t=0$, all the cosines in
\eqref{BN2} have the same phase, adding up constructively yielding
a growth of $B_N = O(N)$.
This will spoil any convergence according to our previous analysis.
When $t \ne 0$, but fixed, when we let $N$ be arbitrarily large,
the same thing will happen since $g_N / p_N \to 0$. Therefore, we expect
the CLT to be valid for $t \ne 0$, but only up to a certain range
$N \le N_c$, such that for $n \le N$, $t g_n/ p_n$ is still big
enough to create difference in the phases to allow cancellations.
We can already see that if $t$ is large enough to compensate
the decaying of $g_n/p_n$, this will be possible. Thus $t$ and
$N_c$ must be related.

The growth of \eqref{BN2}
should be seen even through a smooth approximation.
This growth does not come from the fluctuations in the primes.
Estimating this series through the PNT we have
\beq
B_N(t) = \int_{2}^{p_N} \cos\(t \log x\) d \pi (x) 
    \sim \int_{2}^{p_N} \cos\(t \log x\) \dfrac{dx}{\log x}. \label{BnLog}
\eeq
In the limit $x\to\infty$ we have $\tfrac{\pi(x)}{x/\log x} \to 1$,
so we should interpret our estimate as being asymptotic and for
really large $p_N$. The integral \eqref{BnLog} is easily solved and equal to
\beq
\dfrac{1}{2}\(\Ei\big[\alpha(x)\big] + 
  \overline{\Ei\big[\alpha(x)\big]} \)\Big\vert_{2}^{p_N}
\approx
\dfrac{1}{2} \(\Ei\big[\alpha(p_N)\big] +
 \overline{\Ei\big[\alpha(p_N)\big]}\)
\eeq
where $\alpha(x) = (1-it)\log x$, and $\overline{\Ei(z)}=\Ei(\overline{z})$.
Using the asymptotic expansion
\beq
\Ei(z) \sim \dfrac{e^z}{z}\(1 + \dfrac{1}{z} + O\(\dfrac{1}{z^2}\)\)
\eeq
we finally obtain
\beq\label{BnEstimate}
B_N(t) \sim \dfrac{p_N}{\log p_N} \dfrac{t}{1+t^2} \sin\( t\log p_N \).
\eeq
The above approximation describes
accurately the growth of the series \eqref{BN2}, but only for very
large $N > N_c$. It cannot
describe the series for $N < N_c$, since in this range the series
is governed by fluctuations in the primes, which the
formula \eqref{BnEstimate} does not capture since it is a smooth
asymptotic approximation.

Note that the amplitude of \eqref{BnEstimate} grows with $N$, while
decays with $t$. To avoid this growth, we need to balance these
two variables, which implies a truncation of the series \eqref{BN2} at
a specific cut-off $N_c$. We can impose the condition
$B_N = O(\sqrt{N})$, i.e. $|B_N| = K \sqrt{N}$ for some constant $K>0$,
by constraining the amplitude of \eqref{BnEstimate}. Then we have
\beq
\dfrac{N \log N}{\log\(N\log N\)}\dfrac{1}{|t|} =
\dfrac{N}{\( 1 + \log\log N/\log N \)}\dfrac{1}{|t|}
\approx \dfrac{N}{|t|} < K \sqrt{N}
\eeq
and thus we obtain
\beq\label{cutoff}
N_c \sim t^2.
\eeq
Therefore, for principal characters, 
the upper bound $B_N = O(\sqrt{N})$ is valid only 
in the range $N \le N_c$.
For low $t$, $N_c$ makes short truncations
in the series,
while for large $t$ the cut-off allows us to sum many terms
without crossing the bound. In other words, $N_c$ \emph{regularizes}
our divergent series $B_N$. Note that \eqref{cutoff} is consistent with
Gonek's result \cite{Gonek1}, discussed in the introduction,
which was obtained in a different way. 
Next, we are going to see that the true
origin of this cut-off is the pole at $s=1$.

\subsection{The connection with the pole}

The fundamental origin of the cut-off $N_c$ can be understood
through the analytic continuation of
\eqref{Series1}. Let us consider the specific case of the $\zeta$-function,
since analogous arguments apply to principal Dirichlet $L$-functions.
The analytic continuation of \eqref{Series1}, for $\sigma > 0$,
is given by the formula \cite{ZetaP, ZetaP1}
\beq\label{AZetaP}
X(s) = \sum_{n=1}^{\infty}\dfrac{\mu(n)}{n} \log \zeta(ns).
\eeq
Since $\zeta(s)$ has a pole at $s=1$, when $s=1/n$, $X(s)$ has a
singularity. Note that on the right-half part of the critical strip,
we have only two singularities.
One at $s=1$ and the other at $s=1/2$.
On the left-half part, however, there are an infinite number of
singularities and they accumulate near the point $s=0$.

In the region $1/2<\sigma < 1$, when
$t\to 0$, the cut-off \eqref{cutoff} allow us
to sum only very few terms, avoiding the divergence
for low $t$. Away from the real line, the series for $X(s)$ 
should behave more and more like a strictly convergent series, 
then the cut-off \eqref{cutoff} allows us to sum
a large number of terms.
Therefore, $N_c$ is fundamentally related to the pole of
$\zeta(s)$ at $s=1$,
since it is the origin of these divergences on
the real line segment $0<s<1$.
Away from the pole, for large $t$, the cut-off does not impose 
severe constraints on a truncation of the series \eqref{Series1}.

For non-principal Dirichlet $L$-functions,
the analog of the analytic continuation \eqref{AZetaP} contains
$\log L(ns,\chi)$ instead of $\log \zeta(ns)$.
Since $L(s,\chi)$ has no pole, all the previous discussion does
not apply. In this case there is no need for a cut-off $N_c$, and 
we can take $N\to \infty$ regardless of $t$.

\subsection{Numerical verification}

\begin{figure}[t]
\centering\includegraphics[width=0.5\textwidth]{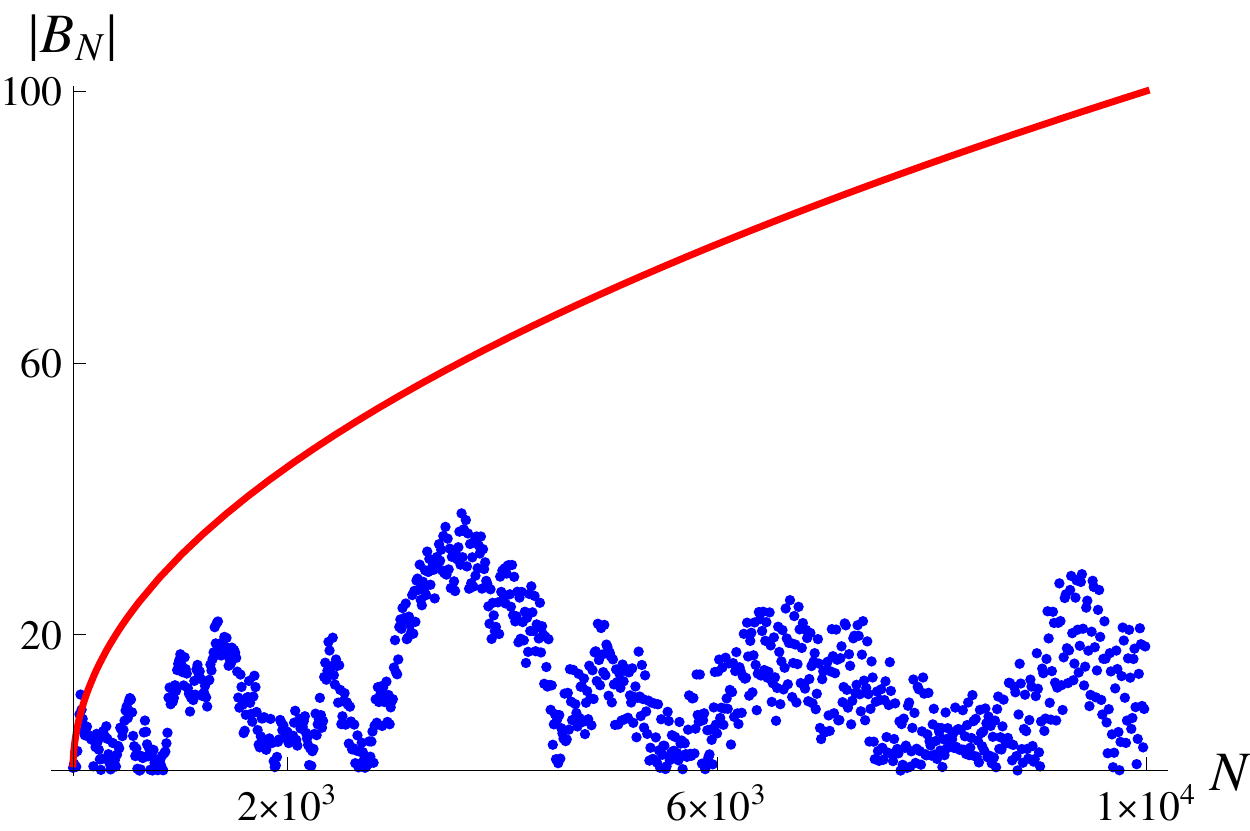} 
\caption{The absolute value of the partial sum 
$B_N = \sum_{n=1}^N \cos\( t \log p_n \)$ 
versus $N$ for a fixed $t=10^3$ (blue dots). The solid red curve 
is $\sqrt{N}$.}
\label{BN_sqrt} 
\end{figure}

The convergence of $B_N / \sqrt{N}$  then does not rely on any 
special detailed properties of the primes,  but rather the opposite, 
on  their multiplicative independence.
Let us numerically test our previous conclusions for \eqref{BN2},
related to the $\zeta$-function.
Obviously, the same is true for non-principal Dirichlet $L$-functions
through \eqref{Bprime}. In this case the numerical 
evidence is even better, and we do not need a cut-off.
Let us check the growth for \eqref{BN2} predicted
by the previous argument.
In Figure~\ref{BN_sqrt}  we plot the partial sums 
$|B_N|$  and one clearly sees this  
$\sqrt{N}$ growth. Note that we choose a range
$N < N_c \sim t^2$, as discussed before.
 
\begin{figure}[t]
\begin{minipage}{0.49\textwidth}
\centering\includegraphics[width=.9\textwidth]{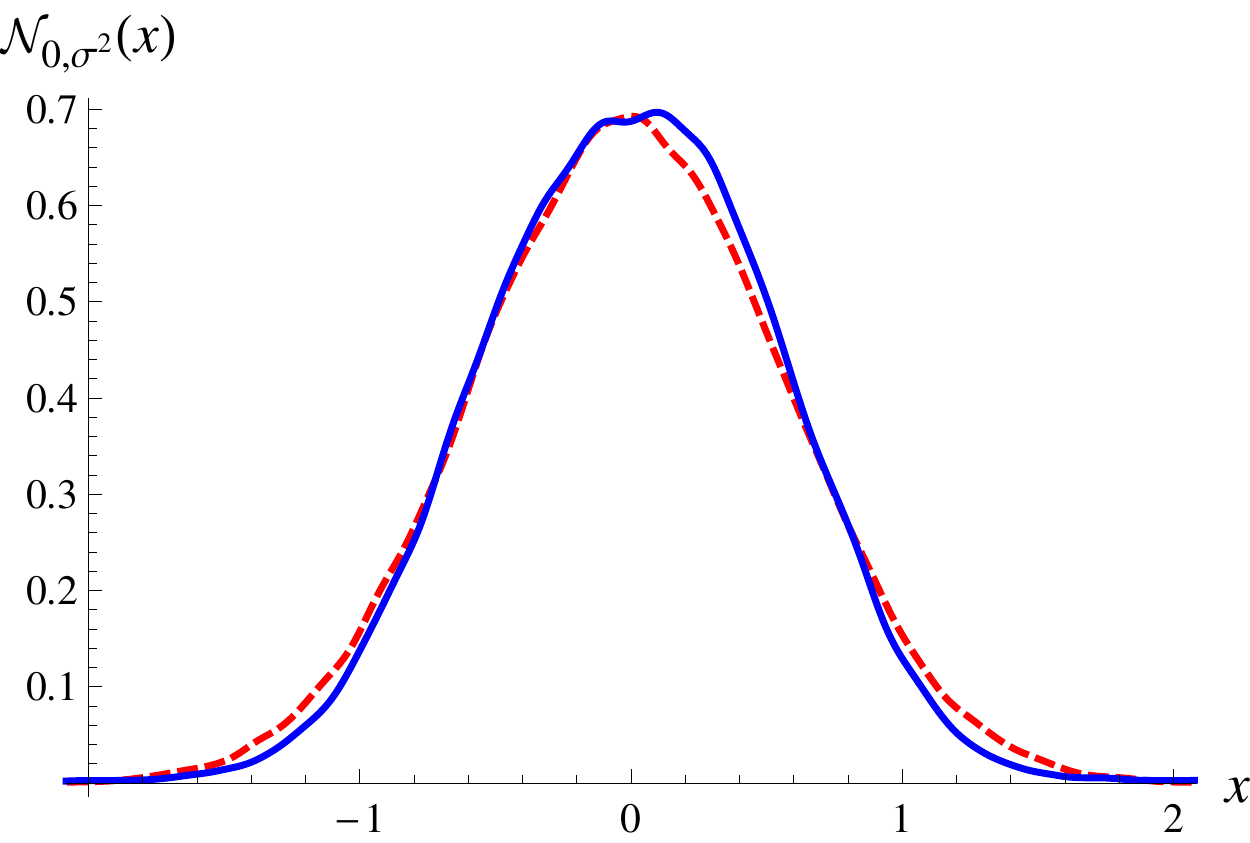}
\end{minipage}
\begin{minipage}{0.49\textwidth}
\centering\includegraphics[width=.9\textwidth]{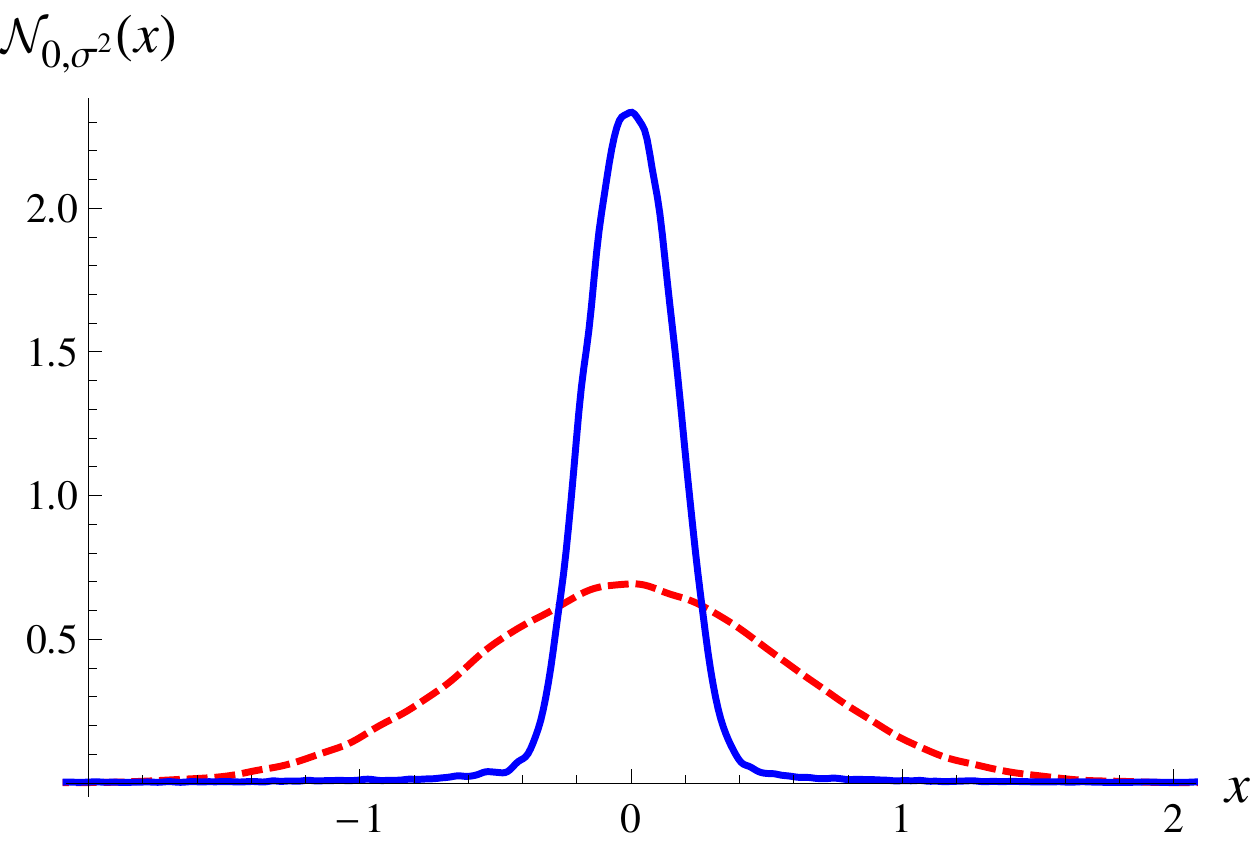}
\end{minipage}
\caption{\textbf{Left:} The probability density function
of $R_N/\sqrt{N}$ (red dashed line) 
compared with the density of $B_N /\sqrt{N}$ (solid blue line). 
We used $t=10^3$, $N=3 \cdot 10^4$ 
and $E=8 \cdot 10^4$ ensembles. Numerically these distributions 
are a gaussian with variance $\sigma^2 \approx 0.578$ 
and  approximately zero mean.
\textbf{Right:} exactly the same numerical experiment
but with the replacement $p_n \approx n \log n $ in \eqref{Bprime}. 
We clearly see that with this approximation 
the CLT is no longer  valid.}
\label{GaussianL} 
\end{figure}

Let us also confirm the gaussian distribution as stated
in \eqref{prob},
using the additional freedom that comes from the random variable $u$.
It is important to note that $u$ is not chosen 
independently for each cosine term  in
the sum; rather one chooses a fixed $u$, randomly, then computes 
the sum $B_N (u)$. In this sense a single sum $B_N(u)$ is completely
deterministic for a given $u$.   
Now consider an ensemble $\{B_N(u_i)/\sqrt{N}\}_{i=1}^{E}$, where for each
element of the set we choose a random $u_i$ in the interval $[0,2\pi]$. 
Now we  can consider its density distribution.
In  Figure~\ref{GaussianL} (left)  we plot 
the this density for both $R_N/\sqrt{N}$ and $B_N/\sqrt{N}$, 
where $R_N$ is given by \eqref{Rseries} with $r_n$ a random variable
on $[-1,1]$. One sees that they are nearly indistinguishable.
This is compelling numerical evidence  that $B_N / \sqrt{N}$ 
approaches a gaussian distribution at large $N$, with zero mean and 
variance $1/2$, since $R_N / \sqrt{N}$ obeys the CLT with a normal
distribution $\tfrac{1}{\sqrt{\pi}}e^{-u^2}$.  
In Figure~\ref{GaussianL} (right) we can
also see that if we use the PNT and replace $p_n \approx n \log n$ 
in \eqref{Bprime},
we loose this normal distribution.  
This is expected since with this approximation
the $\lambda_n$'s will not be linear independent.
This shows that the CLT for \eqref{Bprime}, or \eqref{BN2},
comes from fluctuations in the primes, which are not captured
by a smooth approximation. 

\section{Some consequences of the Euler product formula }
\label{sec:transcendental}

Having provided analytical arguments and numerical evidence 
(see Appendix~\ref{sec:numerical}),  
in this section  we assume the EPF is valid in the  sense described above  
for $\Re (s) > 1/2$, and discuss  some
possible  consequences.     As  already stated in the introduction,   
one  consequence is the validity of the RH,   and this extends to the 
GRH for Dirichlet $L$-functions.       For simplicity,  we limit this
discussion to the $\zeta$-function,  however it easily extends,
and even more precisely,   to the non-principal  Dirichlet $L$-functions,
since no cut-off $N_c$ is required.  

\subsection{The function $S(t)$}
Let $N(T)$ denote the number of zeros  \emph{in the entire critical strip},
$0< \sigma <1$, up to height $T$, where $T$ is not the ordinate 
of a zero.     There is a known exact formula for $N(T)$  due to 
Backlund \cite{Edwards},
\beq
\label{NT}
N(T) = \inv{\pi} \vartheta (T) + 1 + S(T),
\eeq
where 
$\vartheta (T)$ is the 
Riemann-Siegel $\vartheta$ function
\beq
\label{RS}
\vartheta (T) = \arg \Gamma(\tfrac{1}{4} +  
i\,  \tfrac{T}{2})  - T \log \sqrt{\pi}
\eeq
and 
\beq
\label{St}
S(T) =   \inv{\pi}  \arg \zeta \(\tfrac{1}{2} +   i T\).
\eeq
This  result is obtained by the argument principle.
Here,  $S(T)$ is defined by piecewise integration of $\zeta'/\zeta$ from
$s=2$, to $2+i T$,  then to $1/2 + i T$.    
$N(T)$ is a monotonically increasing staircase function,  
however it  is discontinuous 
at the ordinates of non-trivial zeros, where it jumps by the 
multiplicity of the zero.   
Since $\vartheta (T)$ is smooth,
these jumps come from $S(T)$.     

Now,  if the EPF is valid,  then  there are no zeros to the right 
of the critical line.  Then 
$S(T)$ defined by
piecewise integration does not encounter any zeros as one approaches the 
critical line  in the piecewise integration, and  must be the same as
\beq
S(T) = \lim_{\delta \to 0^+}  S_\delta (T)
\eeq
where
\beq
\label{SNew}
S_{\delta}(t) \equiv
\dfrac{1}{\pi}\arg\zeta\(\smallhalf+\delta+it\) = 
-\dfrac{1}{\pi}\lim_{N\to\infty}
\Im\[ \sum_{n=1}^{N} \log\(1 - p_n^{-1/2-\delta-it}\)\].
\eeq
This is an explicit formula for $S(T)$ expressed as a sum over primes,  and    
for  $\delta$ strictly not zero, 
$S_\delta (T)$ is continuous.   As explained in \cite{FrancaLeclair},  
$S(t)$ defined by this limiting procedure is also well-defined at the 
ordinate of a zero on the critical line.   

\begin{figure}[t]
\centering\includegraphics[width=.5\textwidth]{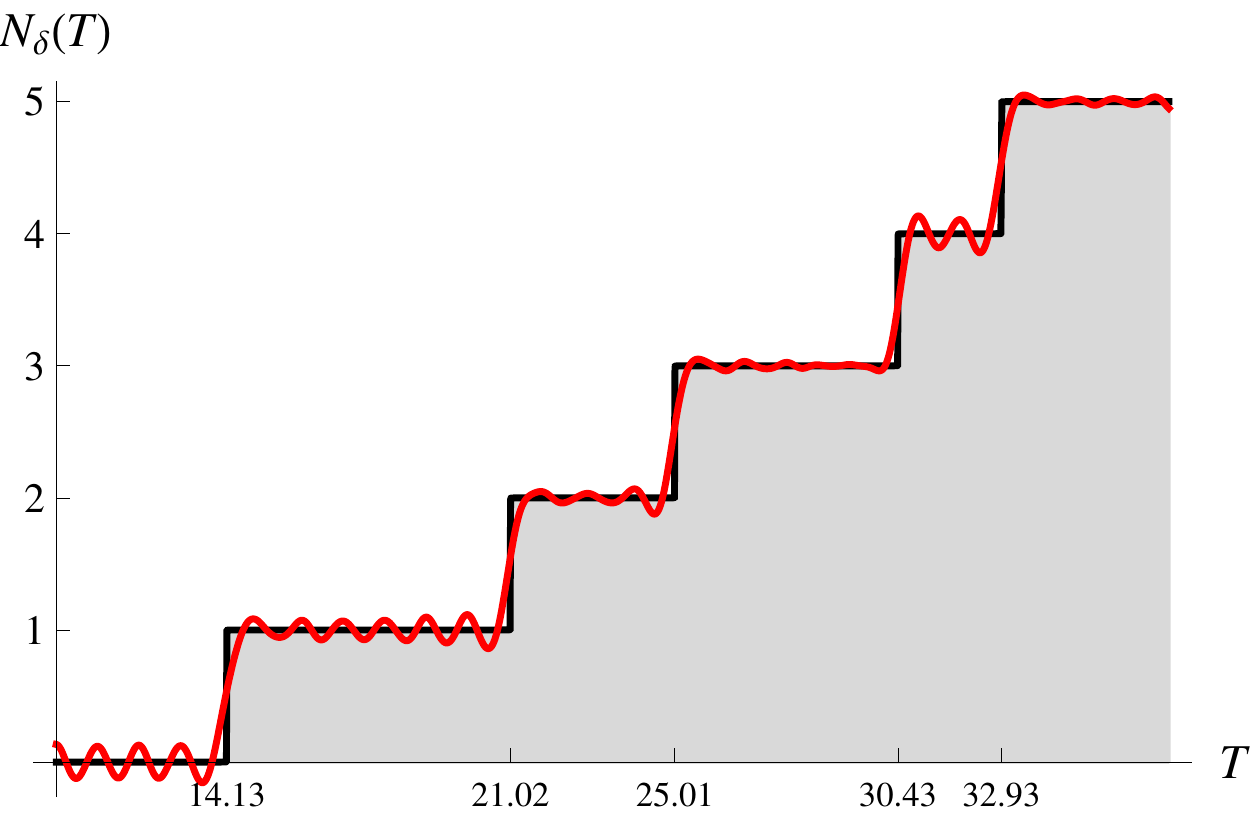}
\caption{The red line shows a plot of equation 
\eqref{NT}  with $S(T)$ replaced by $S_\delta(T)$ defined in terms of 
primes in \eqref{SNew}, i.e.
$N_\delta (T)  \equiv  \inv{\pi} \vartheta(T) + S_\delta(T) +1$.
Here $\delta = 10^{-3}$.  
We use only $100$ primes in the sum. 
The solid black line
corresponds  to the same equation but with the exact 
$S_\delta(T) = \tfrac{1}{\pi} \arg\zeta\(\tfrac{1}{2}+\delta+iT\)$.
Since $N_\delta (T)$ jumps by one at each zero,  
this indicates  that each individual non-trivial zero is related
to an infinite sum over primes.}
\label{zerosprimes}
\end{figure}

The function $S(T)$ knows about the Riemann zeros since it jumps 
at each zero.   Thus,  the expression \eqref{NT}  for $N(T)$,  with $S(T)$ 
replaced by $S_\delta (T)$  in  \eqref{SNew}, which involves a sum over 
primes,  is  a relation between  
Riemann zeros and the primes that is completely the inverse
of Riemann's result for the prime number counting function 
$\pi (x)$ expressed as a sum over non-trivial zeros.
For the latter,  one needs to sum over all zeros
to identify the primes.
Our result is the inverse;  to find the zeros,  one must sum over
all primes. In this sense the distribution of non-trivial zeros
on the critical line is explicitly determined by the prime numbers.
In Figure~\ref{zerosprimes} we plot equation 
\eqref{NT} with $S(T) \to S_\delta(T)$, 
given by \eqref{SNew}, with a finite (small) number of primes. 
The jumps correspond to the nontrivial zeros of $\zeta(s)$.
A stronger version of this is presented below;  see the discussion
following \eqref{transeq}.
If one replaces $p_n\approx n\log n$,  then $S_\delta(T)$ no longer 
jumps at the zeros.
This indicates that the zeros themselves and their GUE statistics
\cite{Montgomery,Odlyzko2}  arises
from fluctuations in the primes.   It would be very interesting to understand
the origin of the GUE statistics in this way.

\subsection{A transcendental equation for the $n$-th zero}
Let us  characterize precisely the zeros on the upper-half of the 
critical line,
$\rho_n = \tfrac{1}{2} + i t_n$ for
$n=1,2,3,\dotsc$. In \cite{RHLeclair,FrancaLeclair}  a transcendental 
equation for 
each $t_n$ was proposed which depends only on $n$. A more lengthy 
discussion of this result can be found in our lectures \cite{FL3}.  
This  transcendental equation for $t_n$ is simple to describe.
We are going to consider only the $\zeta(s)$ case to be concise, 
but the same is easily extended to Dirichlet $L$-functions.
Let $\theta (s) = \arg \chi (s)$ where $\chi$ is the 
completed $\zeta$-function defined in \eqref{chi}.   
It was argued that the  zeros are in one-to-one correspondence 
to the zeros of $\cos \theta$, namely
\beq
\label{thetaeq} 
\lim_{\delta \to 0^+}  \theta \(\tfrac{1}{2} + \delta + i t_n \) =  
\(n-\tfrac{3}{2}\) \pi .
\eeq  
As explained in \cite{FrancaLeclair},   
if the above equation has a unique solution
for every $n$,   then the RH is true and all zeros are simple.     
However,  in that work we  were unable to prove that this equation has a unique solution 
for every $n$.    As we now describe,  the EPF  helps to  resolve these issues.  
Let us first provide a different derivation of \eqref{thetaeq} based on 
the EPF. Using $S_\delta (T) $ in  \eqref{NT},  $N(T)$ is now a monotonically 
increasing staircase function that is smoothed out at the jumps, i.e. it is 
continuous everywhere (see Figure~\ref{zerosprimes}).  
Since it jumps at the ordinate of a zero $t_n$,  and the EPF implies 
there are no zeros off the critical line,  one can use
$N(T)$ to find an equation for $t_n$. Assume  for the moment that 
all zeros 
are simple.  (The derivation of the equation \eqref{thetaeq} in 
\cite{FrancaLeclair}  did not assume this.)   
Then one simply replaces $T \to t_n$ 
and $N\to n-\smallhalf$ in $N(T)$:
\beq
\label{transeq} 
\vartheta \(t_n\) +\pi \lim_{\delta\to0^+}S_\delta\(t_n\) = 
\(n-\tfrac{3}{2}\)\pi .
\eeq
This equation is identical to \eqref{thetaeq}.
The small $\delta$ is required to be positive because the EPF is only valid to
the right of the critical line.
The EPF combined with the properties of $N(T)$  implies 
that the left hand side of the above equation is monotonic and 
continuous,  thus  there is a unique solution to \eqref{transeq}  
for every $n$.

Using the above definition \eqref{SNew}  for $S_\delta$ in terms of 
primes,  the above 
equation \eqref{transeq}
no longer makes 
any reference
to the $\zeta$-function itself.    
This indicates that every single  individual zero depends on
all of the primes.   We were actually able to calculate zeros from 
\eqref{transeq} and \eqref{SNew}.  For instance,  
for the $n=10^5$,  with $N=10^4$ primes,   we obtained
$t_n \approx 74920.826$  whereas the actual value is 
$t_n \approx 74920.827$.  

The  $S_\delta(t)$ term in \eqref{transeq}  fluctuates 
and is very small compared with the 
$\vartheta (t)$ term for large $t$.   If one ignores it, 
and uses  Stirling's approximation  for  
the $\Gamma$-function,  then the solution to the resulting equation can be
expressed in terms of the Lambert $W$-function \cite{FrancaLeclair}:
\beq
\label{Lamb}
t_n \approx  \frac{ 2 \pi (n-\tfrac{11}{8} )}{W[e^{-1} (n-\tfrac{11}{8} )]}.
\eeq
The equation \eqref{transeq}  was used to numerically calculate many zeros
to very high accuracy,  thousands of digits,  
up to the billion-th zero \cite{FrancaLeclair}.   
The  approximation  \eqref{Lamb} is also quite accurate;  
generally the integer part 
is correct,  but it is smooth and does not capture the 
fluctuations that satisfy 
GUE statistics.  
This is clear since this approximation does not capture any sum over primes.    
This suggests that the GUE statistics of the zeros originates from 
the fluctuations 
of the primes.   

As previously stated,  the equation \eqref{transeq} is identical to 
the equation \eqref{thetaeq} 
which comes from $\cos \theta =1$.    
In \cite{RHLeclair,FrancaLeclair}  
the argument
which led to \eqref{transeq}  was  entirely different than the one 
presented here,
i.e. it did not assume the RH nor the simplicity of the zeros,  
and did not rely on  the EPF nor 
knowledge of $N(T)$.   It was obtained directly on the critical line  
using the functional equation. 

The above discussion extends to Dirichlet $L$-functions.   
The analogs  of the above  transcendental 
equations and the Lambert $W$ approximation  for Dirichlet $L$-functions,  
and $L$-functions 
based on modular forms,   were  already presented in \cite{FrancaLeclair}.

\subsection{A Counterexample}

A  well-known counterexample to the RH  is based on the 
Davenport-Heilbronn
function $\CD(s)$,   which is a linear combination of  
two  Dirichlet $L$-functions of modulus $k=5$, i.e.  
$L(s, \chi_{5,2})$ and $L(s, \bar{\chi}_{5,2})$.    
It satisfies a functional equation like \eqref{chi}.
This function is known to have an infinite number of zeros
on the critical line,  but also has non-trivial zeros 
off the  critical line and inside the critical strip.
The Dirichlet $L$-functions each have an Euler product,  
however the sum does not.
The analog of \eqref{transeq}  was studied for this function in 
\cite{FrancaLeclair,FL3}.   
It was found that the analog of $S_\delta (t)$, i.e. 
$S_{\CD}(t) = \tfrac{1}{\pi}\lim_{\delta \to 0^+}
\arg\CD\(\tfrac{1}{2}+\delta+i\,t\)$,  
becomes ill-defined in the vicinity of ordinates $t$ corresponding to
zeros off of the critical line, and there are no solutions to the 
analog of \eqref{transeq} at these points.
This is now perfectly clear,  since there is no Euler product formula 
to smooth out $S_\CD$ here. 

This example provide further evidence that    
the validity of the RH depends on both the \emph{functional equation}
\eqref{chi}  and the \emph{Euler Product Formula} \eqref{Eu}.

\section{Conclusions}
\label{sec:conclusion}

In this paper we provided arguments that the average, or more specifically
the Ces\` aro average, of the Euler product converges into the right-half
part of the critical strip. This implies that the Euler product itself
is meaningful in this region.  Extensive numerical support  for
our statements
are presented in Appendix~\ref{sec:numerical}.

The most important, and delicate,  argument, 
is the central limit theorem applied to the random walk
of the primes, the series  \eqref{BNdef}.
This sum behaves like a normal random variable because of the multiplicative
independence of the primes, and the strongly multiplicative property
of Dirichlet characters. Furthermore, for \emph{non-principal} Dirichlet
characters, the central limit holds as stated, allowing  an
arbitrarily large $N$ limit. However, the situation is more subtle
for  \emph{principal} Dirichlet $L$-functions, 
including the $\zeta(s)$ function.
In these cases, there is no contribution from the character to
\eqref{BNdef}, and because the gaps between primes get smaller
in comparison to the primes, we need to introduce a cut-off, truncating
the series \eqref{cutoff}.
 However, for large
$t$, the cut-off does not impose severe constraints 
since $N_c \sim t^2 \to \infty$.
The cut-off avoids the regime where the cosines sum up
constructively. Note that this regime is not due to fluctuations in
the primes, but rather the opposite. 
Furthermore, the cut-off is intimately related to the pole
at $s=1$. Our analysis in effect predicts the existence
of a pole.  Since non-principal Dirichlet $L$-functions
have  no such a pole, there is no need for a cut-off in this case.

We also discussed some consequences assuming that the Euler product
formula is valid for $\Re(s) > 1/2$. The most important one, is that
the function $S(t) = \tfrac{1}{\pi}\arg\zeta(1/2+\delta+it)$ becomes smooth
through the $\delta \to 0^{+}$ limit, and thus we expect the transcendental
equation \eqref{transeq} to have a unique solution for
every $n$ \cite{FrancaLeclair}.
The same argument applies to Dirichlet $L$-functions, with their 
respective transcendental equations, also presented in \cite{FrancaLeclair}.

\begin{acknowledgments}
We thank Denis Bernard,  
Keith Conrad,  and Jon Keating for discussions.
We also wish to thank Merlin Enault for 
the last 2 entries in Table \ref{table}.
GF thanks the support from CNPq-Brazil.
\end{acknowledgments}

\bigskip

\appendix

\section{Numerical studies}
\label{sec:numerical}

We have provided arguments that the logarithm of the 
Euler product 
Ces\`aro-converges in  the region $\sigma > 1/2$,  
subject to the qualification for principle characters described in
the last section.  
We now present compelling  numerical evidence of  the validity of this result.  
Throughout this section we plot the Euler product itself,  rather than 
its  average,   since the resolution of the plots 
is not high  enough to see the small  fluctuations,  so that these plots are 
indistinguishable from the plots of the average.  
In other words,  we provide evidence for the following equality
\beq\label{prodpartial}
L(s,\chi) = \lim_{N\to\infty} \langle \CP_N(s, \chi)\rangle, \qquad
\CP_N(s, \chi) = \prod_{n=1}^N\(1-\dfrac{\chi (p_n) }{p_n^{\,s}}\)^{-1},
\eeq
where  $\langle  \CP_N (s,\chi) \rangle$ is its 
arithmetic average over $N$:
\beq
\label{CPAve}
\langle  \CP_N (s,\chi) \rangle  = \dfrac{1}{N}  \sum_{n=1}^N  \CP_n (s,\chi).
\eeq
The above arithmetic average of the product should converge 
for the same reasons that the arithmetic average of $\log \CP(s,\chi)$ 
converges, since averaging is just a smoothing procedure.

\subsection{Riemann $\zeta$-function} 

\begin{figure}[b]
\begin{minipage}{.49\textwidth}
\centering\includegraphics[width=.9\textwidth]{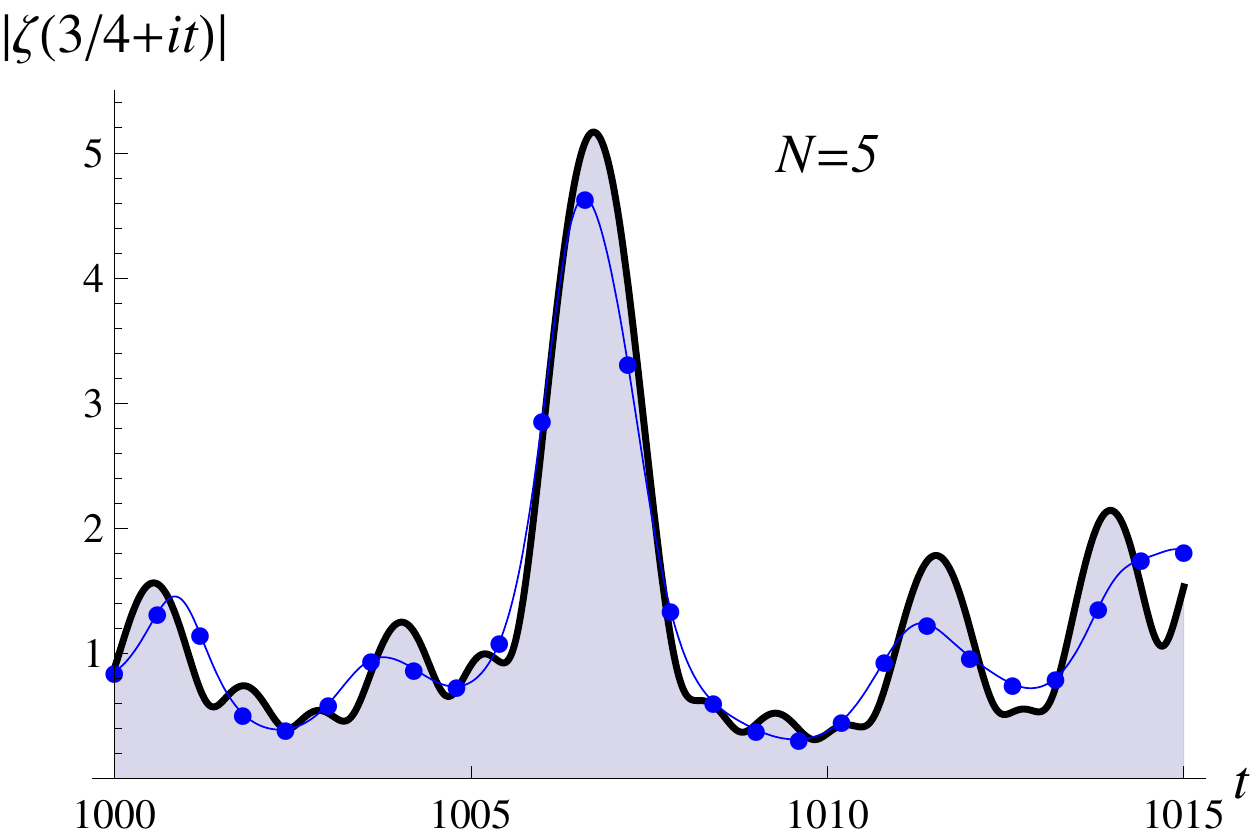}
\end{minipage}
\begin{minipage}{.49\textwidth}
\centering\includegraphics[width=.9\textwidth]{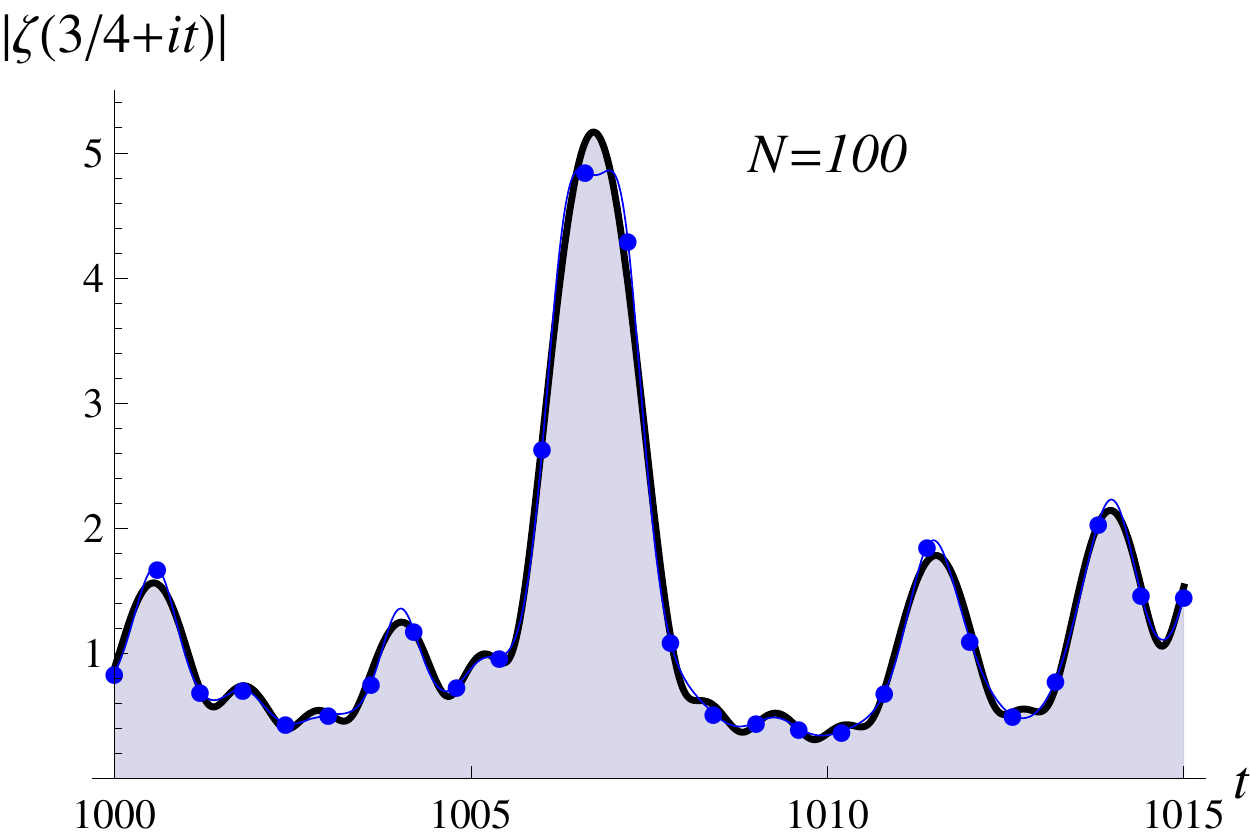}
\end{minipage}
\begin{minipage}{.49\textwidth}
\centering\includegraphics[width=.9\textwidth]{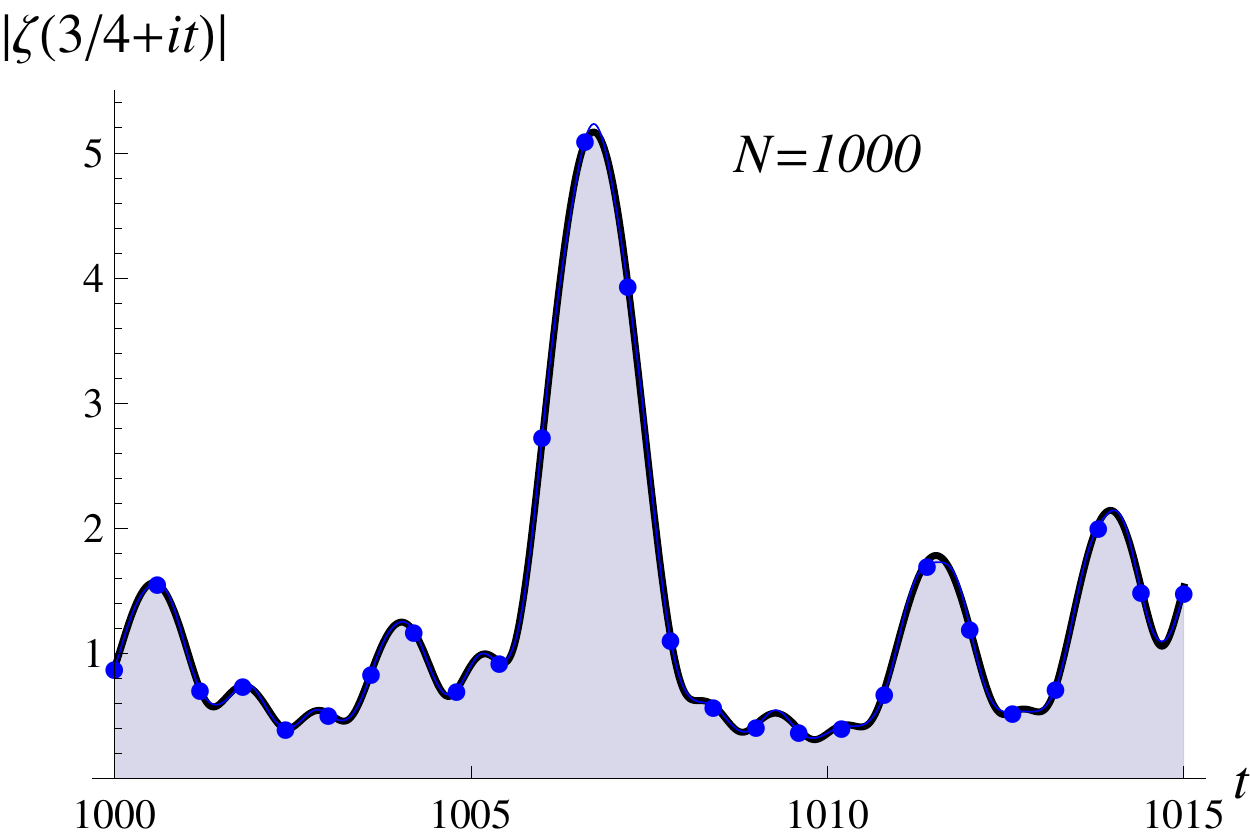}
\end{minipage}
\begin{minipage}{.49\textwidth}
\centering\includegraphics[width=.9\textwidth]{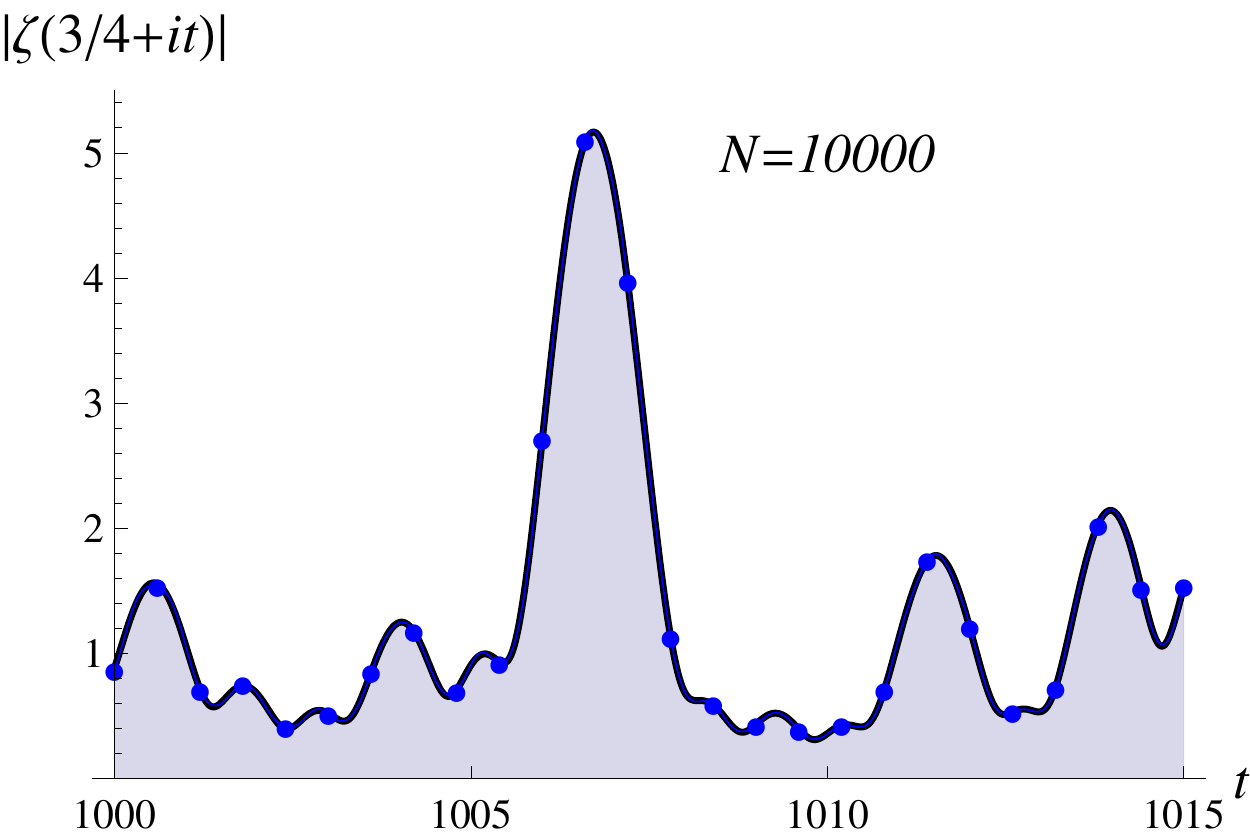}
\end{minipage}
\caption{The black line is the actual  $|\zeta(3/4+it)|$, analytically 
continued into the strip, and
the blue line is the partial product $|\CP_N(3/4+it)|$.
Dots are added to the line to aid visualization.}
\label{abs_zeta}
\end{figure}

\begin{figure}[t]
\centering
\includegraphics[width=0.5\textwidth]{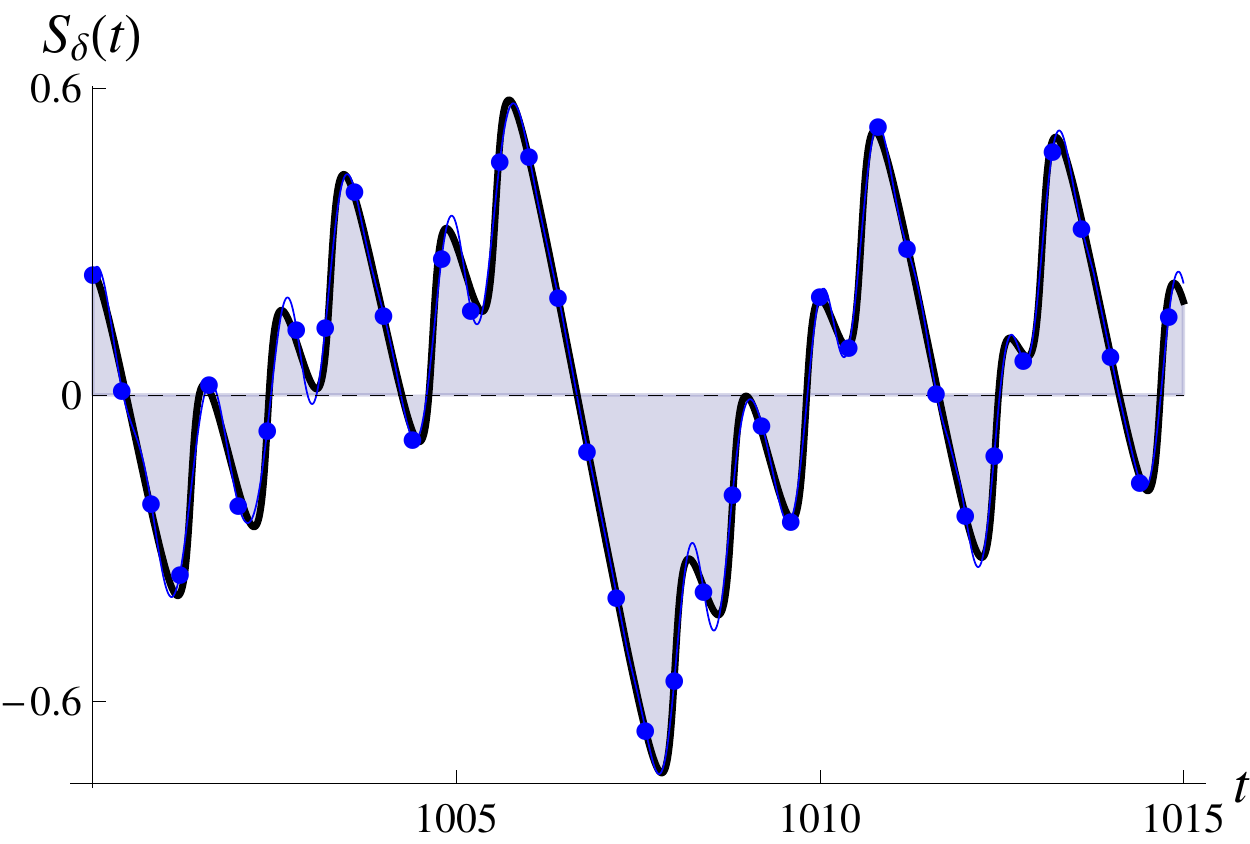}
\caption{The  black line is 
the actual $\tfrac{1}{\pi}\arg\zeta\(\tfrac{1}{2}+\delta+it\)$, 
and the blue line 
is the RHS of \eqref{SNew}. Dots are added 
to the line to aid in visualization.  
We used $\delta=10^{-1}$ and $N=10^5$.}
\label{arg_zeta}
\end{figure}

In Figure~\ref{abs_zeta}
one can see how the partial product in \eqref{prodpartial}
converges to the $\zeta(s)$ function as we increase $N$.
For higher $N$ the curves are indistinguishable. Note, however,
that we cannot go beyond the cut-off \eqref{cutoff} for a given $t$.
Let us also verify convergence for $\arg \zeta$,
which plays a central role for the zeros on the critical 
line (see the Section~\ref{sec:transcendental}).   
Using the EPF we have equation \eqref{SNew}, whose 
equality is verified in Figure~\ref{arg_zeta}.
This assures that both 
the real and imaginary parts of the 
Euler product converge.  
As we approach the critical line $\sigma \to 1/2^+$
higher $N$ is of course  required. 

\begin{figure}[t]
\begin{minipage}{0.49\textwidth}
\centering\includegraphics[width=.9\textwidth]{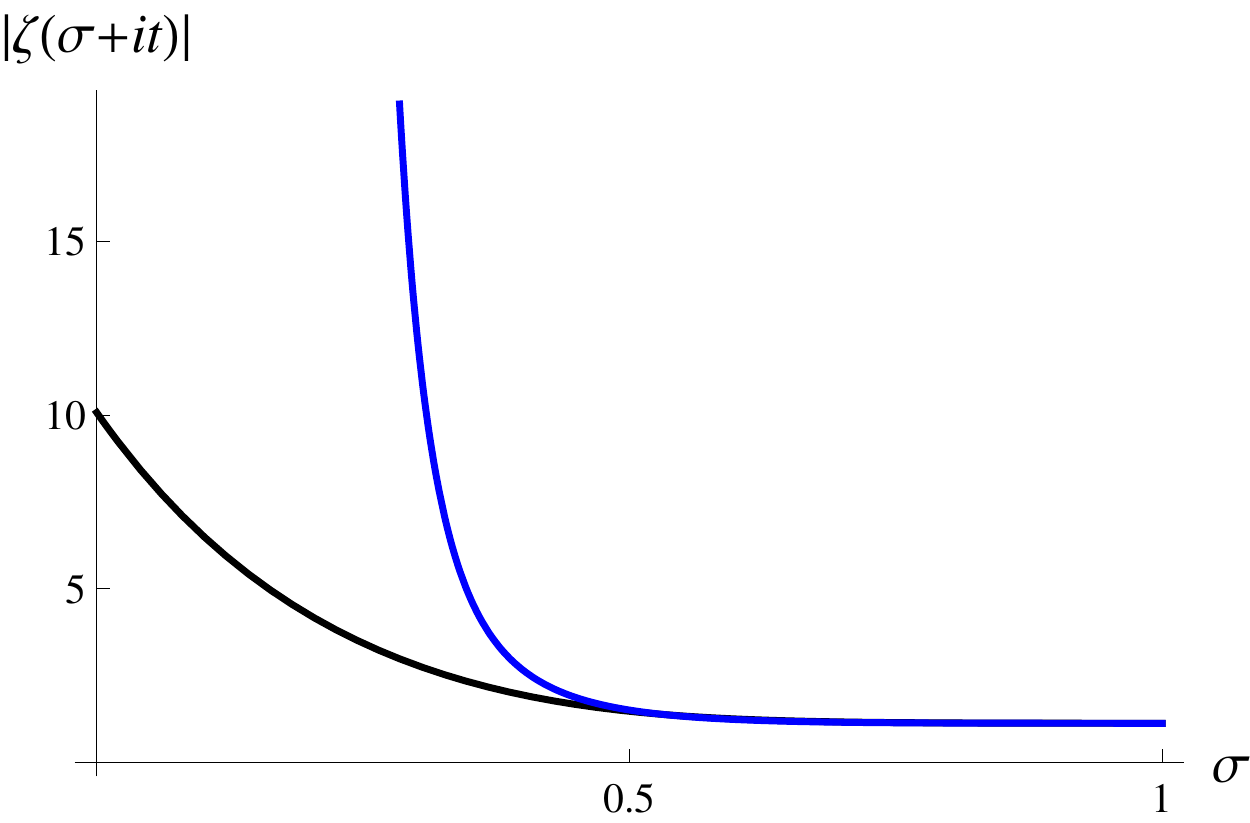}
\end{minipage}
\begin{minipage}{0.49\textwidth}
\centering\includegraphics[width=.9\textwidth]{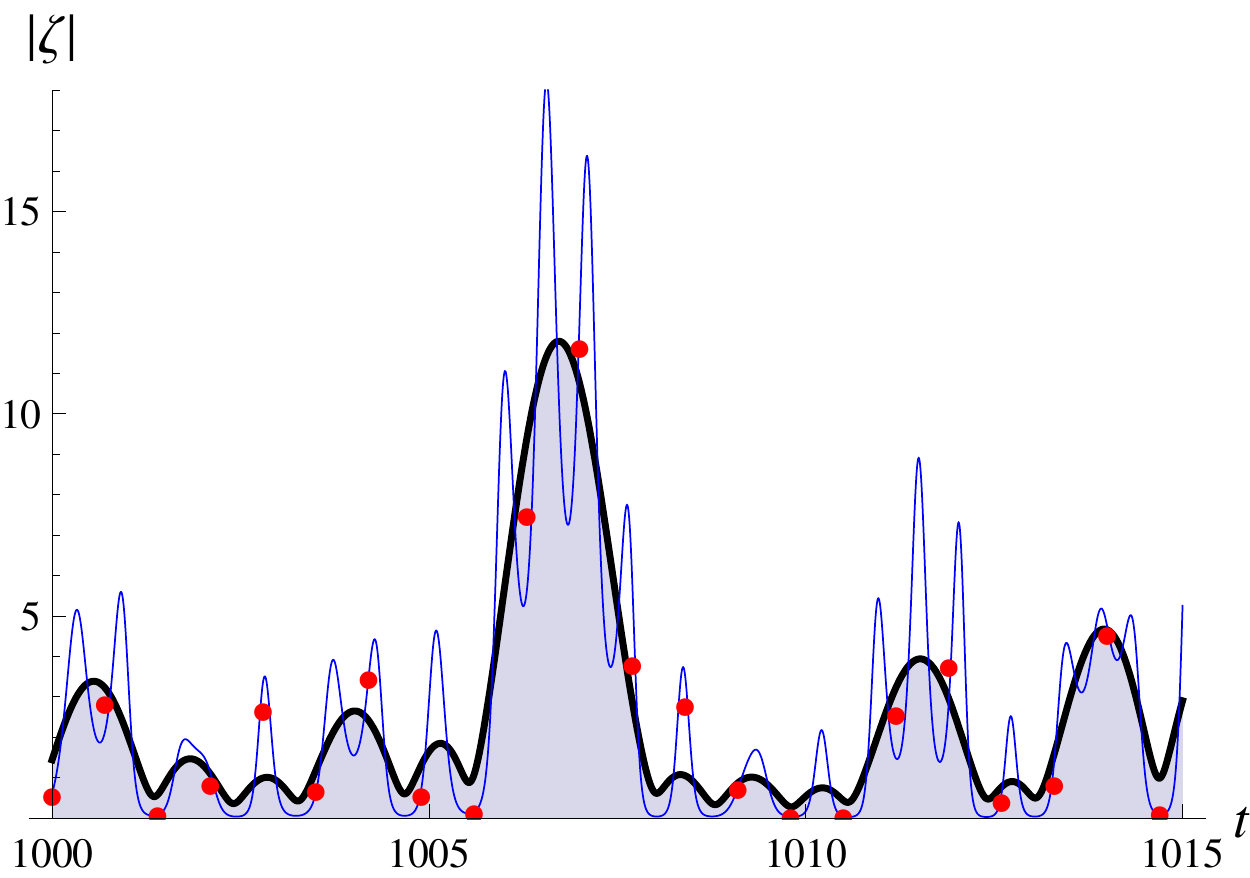}
\end{minipage}
\caption{\textbf{Left:} the black line corresponds to 
$|\zeta(\sigma + it)|$ against $0<\sigma<1$, for $t=500$.
The blue line is the partial 
product $|\CP_N(\sigma + it)|$ with $N=10^4$.
\textbf{Right:} the black line is the exact $|\zeta|$, and the
blue line is the partial product $|\CP_N|$ (with $N=8 \cdot 10^3$), 
against $t$.
We took $\sigma = 0.4$.   The red dots are the Ces\` aro average
$|\langle \CP_N \rangle|$. If we increase $N$ the results are even worse.}
\label{phase_trans}
\end{figure}

One can clearly see how the Euler product formula is not valid for 
$\sigma  \le 1/2$ from Figure~\ref{phase_trans}.
The curves only match for $\sigma>1/2$ and the dramatic change 
in behavior is abrupt at $\sigma = 1/2$,  as predicted. 
The divergences shown in Figure~\ref{phase_trans} (right) get 
worse for higher $N$.

\begin{figure}[t]
\begin{minipage}{0.49\textwidth}
\centering\includegraphics[width=.9\textwidth]{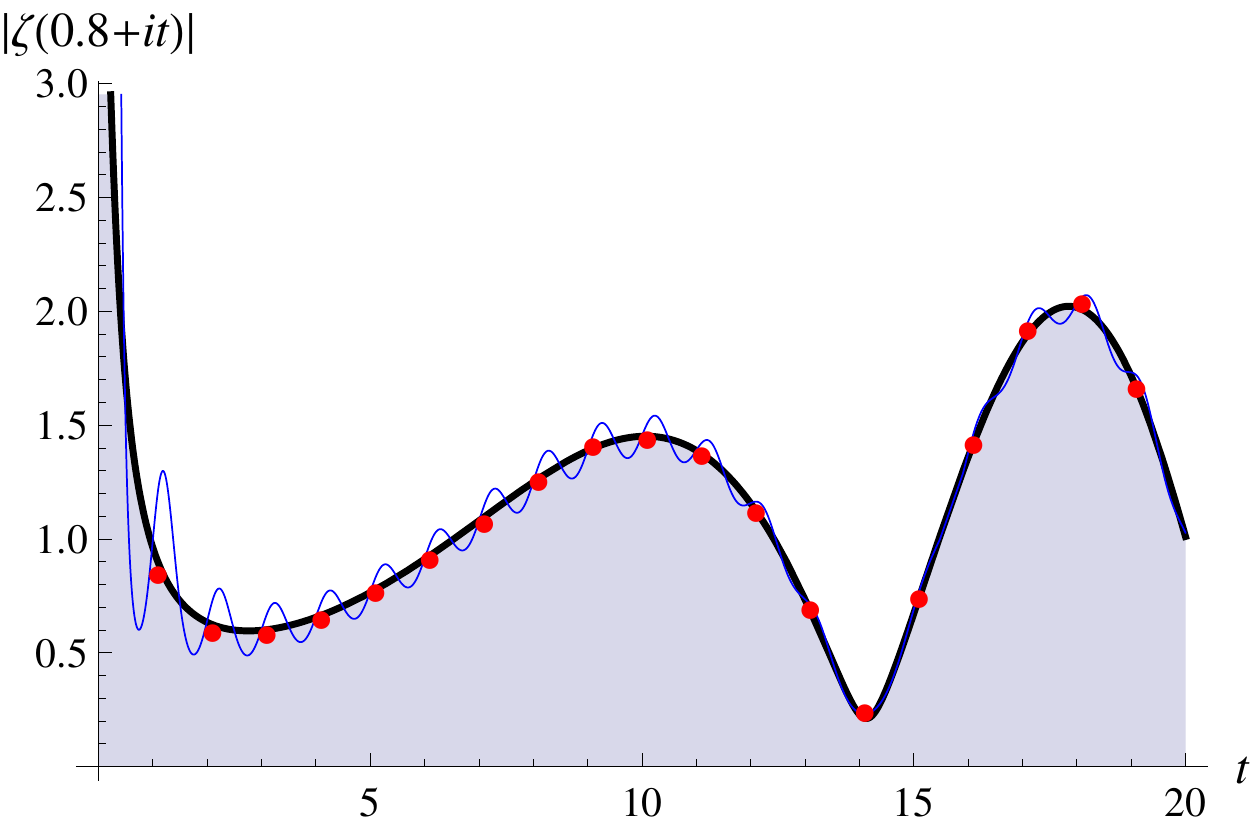}
\end{minipage}
\begin{minipage}{0.49\textwidth}
\centering\includegraphics[width=.9\textwidth]{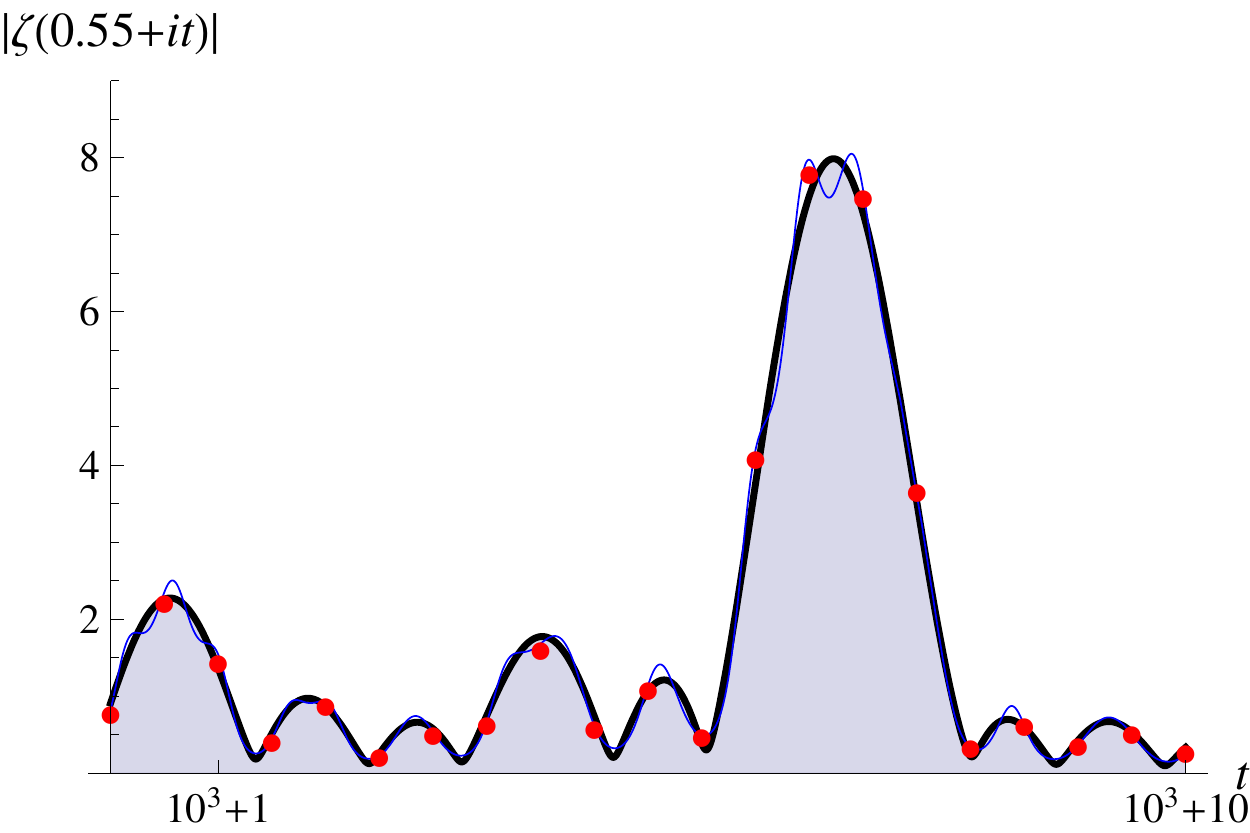}
\end{minipage}
\caption{\textbf{Left:} the black line is 
$|\zeta (0.8 + i t)|$, and the blue line
is $|\CP_N|$, with $N=10^2$,  and the red dots 
correspond to the Ces\`aro average $|\langle\CP_N\rangle|$.
\textbf{Right:} here we have $\zeta(0.55+it)$, and we used $N=4\cdot10^4$. The
convergence close to the critical line requires very large $N$.
}
\label{Lowt}
\end{figure}

As we discussed before, there is no convergence on the real line
$t=0$ due to the pole at $s=1$. 
It is exactly because of this divergence that we had to
introduce the cut-off \eqref{cutoff}. However, for short truncations
of the product, i.e. not so high $N$, 
we can describe the $\zeta$-function quite accurately even
for low $t$, as shown in Figure~\ref{Lowt} (left).
We can see that the oscillations get stronger close to $t=0$, but the
Ces\` aro average is still well-behaved and closer to the actual
value of $\zeta(s)$ than the partial product itself.
The convergence close to the critical line is very slow, and requires
very high $N$. Thus to test the results close to the critical line,
we also have to choose high $t$ due to the cut-off relation \eqref{cutoff}.
One can see from Figure~\ref{Lowt} (right) that the 
Ces\` aro average approximates $\zeta(s)$ correctly,
even close to the critical line.

In Table~\ref{table} we show some values of the average 
$|\langle \CP_N \rangle|$ and the product  $|\CP_N|$ itself. 
The convergence is slow, but
one can see that $\langle \CP_N(s) \rangle \approx \zeta(s)$ as we
increase $N$,  whereas
the unaveraged $\CP_N(s)$ continues to oscillate around $\zeta(s)$. 
With $N=10^5$ we obtain nearly $5$ digit accuracy for $t=100$. 
Note that the results eventually start to get worse
for very high $N$,  here roughly $10^8$.   
We are increasing $N$ much beyond 
the cut-off predicted by \eqref{cutoff}. We can do this since
we are not close to the critical line.

\begin{table}[h]
\def\arraystretch{0.7}
\centering
\begin{minipage}{0.5\textwidth}
\begin{tabular}{@{}ccc@{}}
\toprule[0.8pt]
$N$ & $|\langle \CP_N \rangle|$ & $|\CP_N|$ \\
\midrule[0.4pt]
$1\cdot10^3$ & $0.976752$ & $0.972210$ \\
$2\cdot10^3$ & $0.976690$ & $0.981506$ \\
$3\cdot10^3$ & $0.977653$ & $0.976654$ \\
$4\cdot10^3$ & $0.977865$ & $0.975735$ \\
$5\cdot10^3$ & $0.977926$ & $0.984674$ \\
$6\cdot10^3$ & $0.977463$ & $0.977893$ \\
$7\cdot10^3$ & $0.978208$ & $0.976510$ \\
$8\cdot10^3$ & $0.977593$ & $0.978773$ \\
$9\cdot10^3$ & $0.978290$ & $0.981781$ \\
$1\cdot10^4$ & $0.977900$ & $0.971017$ \\
$1\cdot10^5$ & $0.977703$ & $0.971203$ \\
$1\cdot10^6$ & $0.977925$ & $0.971491$ \\
$1\cdot10^7$ & $0.978168$ & $0.978027$ \\
$1\cdot10^8$ & $0.977823$ & $0.984481$ \\
\rowcolor{gray!25} 
$2\cdot10^8$ & $0.956304$ & $0.885545$  \\
\rowcolor{gray!25} 
$3\cdot10^8$ & $0.924928$ & $0.794254$  \\ 
\midrule[0.4pt]
\multicolumn{3}{@{}c@{}}{$|\zeta(0.95+i\,20)|=0.977848$}\\
\bottomrule[0.8pt]
\end{tabular}
\end{minipage}
\begin{minipage}{0.49\textwidth}
\begin{tabular}{@{}ccc@{}}
\toprule[0.8pt]
$N$ & $|\langle \CP_N \rangle|$ & $|\CP_N|$ \\
\midrule[0.4pt]
$1\cdot10^3$ & $1.690988$ & $1.694894$ \\
$2\cdot10^3$ & $1.692350$ & $1.694156$ \\
$3\cdot10^3$ & $1.692590$ & $1.690354$ \\
$4\cdot10^3$ & $1.692399$ & $1.688480$ \\
$5\cdot10^3$ & $1.691996$ & $1.687150$ \\
$6\cdot10^3$ & $1.691666$ & $1.689158$ \\
$7\cdot10^3$ & $1.691508$ & $1.688145$ \\
$8\cdot10^3$ & $1.691400$ & $1.691700$ \\
$9\cdot10^3$ & $1.691381$ & $1.692973$ \\
$1\cdot10^4$ & $1.691345$ & $1.690480$ \\
$1\cdot10^5$ & $1.691373$ & $1.692136$ \\
$1\cdot10^6$ & $1.691429$ & $1.691577$ \\
$1\cdot10^7$ & $1.691414$ & $1.691703$ \\
$1\cdot10^8$ & $1.691385$ & $1.693287$ \\
\rowcolor{gray!25} 
$2\cdot10^8$ & $1.745257$ & $1.923738$ \\
\rowcolor{gray!25} 
$3\cdot10^8$ & $1.852499$ & $2.203470$ \\
\midrule[0.4pt]
\multicolumn{3}{@{}c@{}}{$|\zeta(0.95+i\,100)|=1.691397$}\\
\bottomrule[0.8pt]
\end{tabular}
\end{minipage}
\caption{
Convergence $\langle\CP_N \rangle$, and 
$\CP_N$, for the $\zeta$-function. Note that even for $N\gg N_c \sim t^2$ 
the results are good, but eventually it starts to
deviate from the correct value as shown in the two last entries.}
\label{table}
\end{table}

\subsection{Non-principle Dirichlet $L$-functions} 

\begin{figure}[b]
\begin{center}
\begin{minipage}{0.49\textwidth}
\includegraphics[width=.9\textwidth]{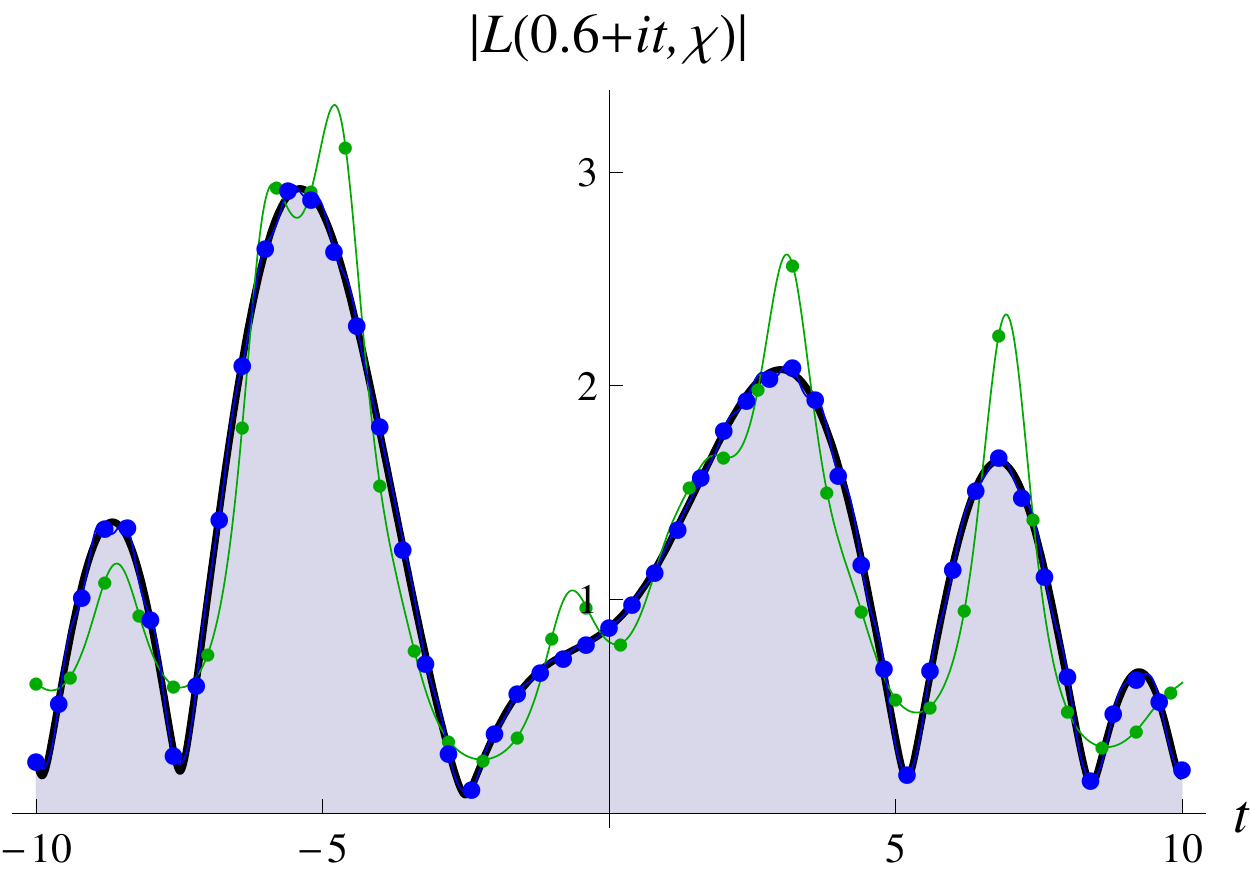}
\end{minipage}
\begin{minipage}{0.49\textwidth}
\includegraphics[width=.9\textwidth]{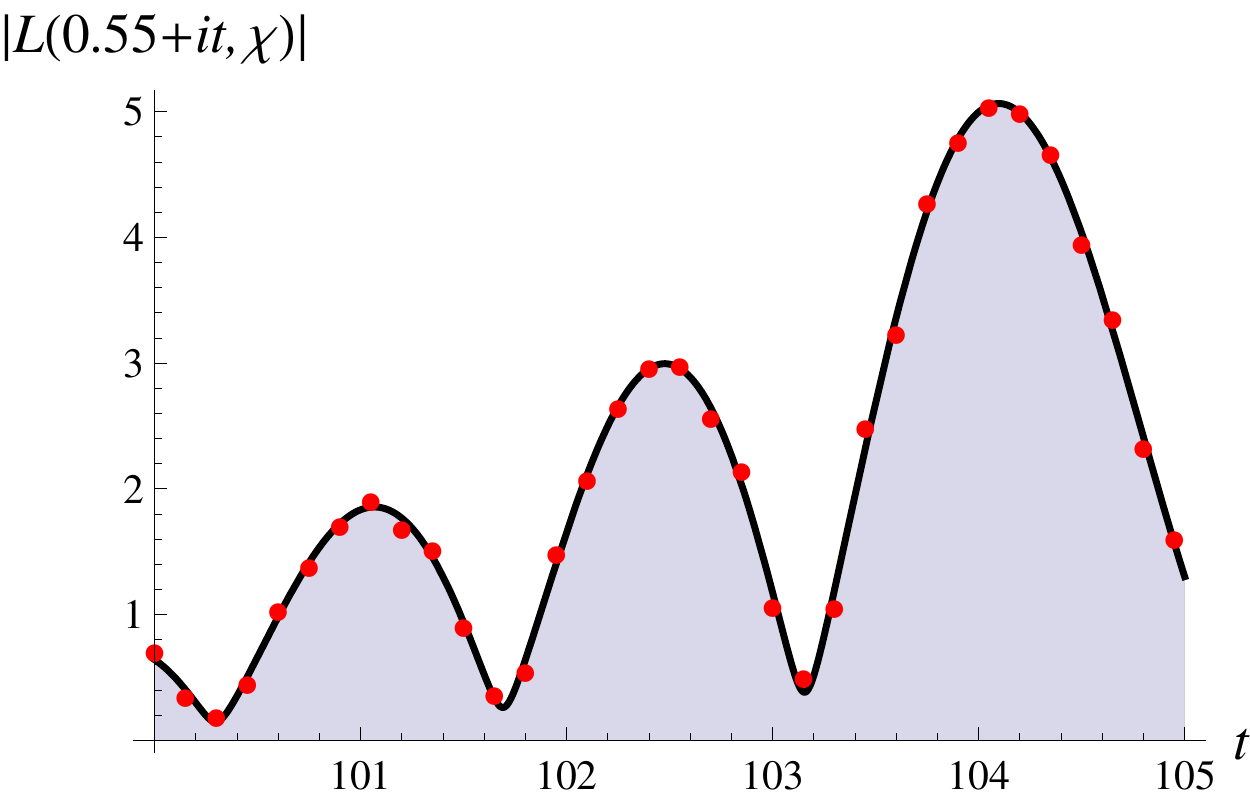}
\end{minipage}
\end{center}
\caption{\textbf{Left:} we have $s=0.6+it$. The black line 
is $|L(s,\chi)|$ with $\chi=\chi_{7,2}$ given in  \eqref{char72}.
The green line is the partial product $|\CP_N(s,\chi)|$ obtained 
with only $N=5$, while the blue line with
$N=10^4$. We can see that the results get better as 
we increase $N$. We can also see that there are no divergences 
on  the real line $t =0$, in contrast with Figure~\ref{Lowt}.
\textbf{Right:} here we are closer to the critical line, $s=0.55+it$.
The black line is $|L(s,\chi)|$, and the red dots correspond to
the Ces\` aro average $|\langle \CP_N(s,\chi) \rangle|$
with $N = 5\cdot 10^6$.}
\label{LowtDir}
\end{figure}

Let us consider a concrete example 
with the primitive character of modulus $7$, shown below:
\beq\label{char72}
\def\arraystretch{1.2}
\begin{tabular}{@{}c|ccccccc@{}}
$n$             & $1$ & $2$ & $3$ & $4$ & $5$ & $6$ & $7$ \\
\hline
$\chi_{7,2}(n)$ &
$1$ & $e^{ 2\pi i /3}$ & $e^{\pi i / 3}$ & $e^{-2\pi i / 3}$ &
$e^{-\pi i / 3}$ & $-1$ & $0$
\end{tabular} 
\eeq
In Figure~\ref{LowtDir} (left) we plot the absolute value of the 
partial Euler product for $L(s,  \chi_{7,2})$, and we can see
how it fits the $L$-function on the right-half part of the 
critical strip, even at
the real line $t=0$, since there is no pole. This is in clear contrast
with Figure~\ref{Lowt} where the Euler product for $\zeta(s)$ is not 
even finite on the line segment $1/2<s\le1$.
Thus one clearly sees that $\zeta(s)$ is
exceptional,  along with the other principal Dirichlet $L$-functions.
The convergence of the Euler product for
non-principal Dirichlet $L$-functions is much better behaved. In fact one
may check that the analog of Figure~\ref{BN_sqrt}  has smaller fluctuations
and resembles  even more  the standard random walk. As explained,
in this case there is no cut-off $N_c$, and we can take $N$ as large as
desirable. In Figure~\ref{LowtDir} (right) we see the Ces\` aro
average in comparison with the actual value of $L(s,\chi_{7,2})$, closer
to the critical line.

In Table~\ref{dirtable} we compute the partial product, and its Ces\`aro
average, for 
$L(s,\chi_{7,2})$. We can see how the last digits fluctuate
but the numbers are close to the actual value of $|L(s,\chi_{7,2})|$.
Note that, contrary to the $\zeta(s)$ case where we are 
limited in accuracy
by the cut-off $N\le N_c$, for non-principal Dirichlet $L$-functions 
we expect the results to get better and better as we arbitrarily increase
$N$.

\begin{table}[h]
\def\arraystretch{0.7}
\centering
\begin{minipage}{0.5\textwidth}
\begin{tabular}{@{}ccc@{}}
\toprule[0.8pt]
$N$ & $|\langle \CP_N \rangle|$ & $|\CP_N|$ \\
\midrule[0.4pt]
$1\cdot10^3$ & $0.8940791$ & $0.8949042$ \\
$2\cdot10^3$ & $0.8947639$ & $0.8951913$ \\
$3\cdot10^3$ & $0.8948319$ & $0.8946522$ \\
$4\cdot10^3$ & $0.8947869$ & $0.8950135$ \\
$5\cdot10^3$ & $0.8948144$ & $0.8946950$ \\
$6\cdot10^3$ & $0.8947834$ & $0.8945271$ \\
$7\cdot10^3$ & $0.8947674$ & $0.8948700$ \\
$8\cdot10^3$ & $0.8947783$ & $0.8947044$ \\
$9\cdot10^3$ & $0.8947768$ & $0.8948476$ \\
$1\cdot10^4$ & $0.8947921$ & $0.8950163$ \\
$1\cdot10^5$ & $0.8949043$ & $0.8949518$ \\
\midrule[0.4pt]
\multicolumn{3}{@{}c@{}}{$|L\(0.95, \chi\)|=0.89492570$}\\
\bottomrule[0.8pt]
\end{tabular}
\end{minipage}
\begin{minipage}{0.49\textwidth}
\begin{tabular}{@{}ccc@{}}
\toprule[0.8pt]
$N$ & $|\langle \CP_N \rangle|$ & $|\CP_N|$ \\
\midrule[0.4pt]
$1\cdot10^3$ & $0.6183514$ & $0.6208759$ \\
$2\cdot10^3$ & $0.6195137$ & $0.6202016$ \\
$3\cdot10^3$ & $0.6199206$ & $0.6211404$ \\
$4\cdot10^3$ & $0.6201229$ & $0.6205615$ \\
$5\cdot10^3$ & $0.6202306$ & $0.6207769$ \\
$6\cdot10^3$ & $0.6202884$ & $0.6205089$ \\
$7\cdot10^3$ & $0.6203365$ & $0.6207366$ \\
$8\cdot10^3$ & $0.6203860$ & $0.6207027$ \\
$9\cdot10^3$ & $0.6204248$ & $0.6207634$ \\
$1\cdot10^4$ & $0.6204524$ & $0.6207338$ \\
$1\cdot10^5$ & $0.6207878$ & $0.6209509$ \\
\midrule[0.4pt]
\multicolumn{3}{@{}c@{}}{$|L\(0.95+i\,100, \chi\)|=0.62101132$}\\
\bottomrule[0.8pt]
\end{tabular}
\end{minipage}
\caption{
Numerical results for the Euler product and its Ces\` aro average
for $\chi = \chi_{7,2}$, shown in \eqref{char72}.
We expect the results to be more accurate as we increase $N$.
We choose $t=0$ in the left table to explicitly show that there are no
divergences on the real line.}
\label{dirtable}
\end{table}

\end{document}